\documentclass[]{article}

\usepackage{stmaryrd}

\usepackage{amsmath,amscd,amssymb}

\usepackage{verbatim}

\usepackage[pdftex, pdfstartview=FitH]{hyperref}

\hypersetup{colorlinks,
        linkcolor=red,
        filecolor=black,
        urlcolor=blue,
        citecolor=green,
        pdftitle={For Alstom},
        pdfauthor={             },
        pdfsubject={            },
        pdfkeywords={           },
        pdfproducer={ps2pdf} }

\newtheorem {theorem}{Theorem}[section]

\newtheorem {prop}[theorem]{Proposition}

\newtheorem {corollary}[theorem]{Corollary}

\newtheorem {lemma}[theorem]{Lemma}

\newtheorem {definition}[theorem]{Definition}

\newtheorem {example}[theorem]{Example}

\newtheorem {remark}[theorem]{Remark}

\newtheorem {notation}[theorem]{Notation}

\def\proof{\par\noindent{\bf Proof.}}
\def\qed{\hfill $\Box$}

\title{Zeta functions of trinomial curves and maximal curves}

\author{Menglong Nie
\\ \\
{\small School of Mathematics, Northwest University,}\\{\small Xi'an,
710127, P.R. China}
\\ \\
{\small\small niemelo@gmail.com}}

\date{}

\begin{document}
\maketitle
\begin{abstract}
We determine the zeta functions of trinomial curves in terms of Gauss sums and Jacobi sums, and we obtain an explicit formula of the genus of a trinomial curve over a finite field, then we study the conditions for a trinomial curve to be a maximal curve over a finite field.
\end{abstract}


\section {Introduction}
\quad Let $C$ be a projective, non-singular, geometrically irreducible algebraic curve defined over the finite field $\mathbb{F}_{q}$ with $q$ elements, let $N_i$ be the number of $\mathbb{F}_{q^i}-$ rational points of $C$, i.e. $N_i=\#C(\mathbb{F}_{q^i})$, then the zeta function $Z(C/\mathbb{F}_{q};t)$ of the curve $C$ is defined as the formal series
\[
Z(t)=Z(C/\mathbb{F}_{q};t):=\exp(\sum_{i=1}^{\infty}\frac{N_i}{i}t^i).
\]
Andr$\acute{e}$ Weil \cite{Weil1} proved that $Z(C/\mathbb{F}_{q};t)=\frac{L(t)}{(1-t)(1-qt)}$ with $L(t)\in \mathbb{Z}[t]$ and $\deg L(t)=2g$, where $g$ is the genus of $C$. The numerator $L(t)$ of the zeta function $Z(t)$ is also called the L-polynomial of $C/\mathbb{F}_{q}$. Further results and information about zeta functions over finite fields can be found, e.g., in \cite{Wan}.

The Hasse-Weil bound for the number of rational points on the curve $C$ of genus $g$ over the finite field $\mathbb{F}_{q}$ states that $\#C(\mathbb{F}_{q})\le q+1+2g\sqrt{q}$. If a curve of genus $g$ over $\mathbb{F}_{q^2}$ attains this bound, then the numerator of its zeta function over $\mathbb{F}_{q^2}$ equals $(1+qt)^{2g}$ and the curve is called a maximal curve over $\mathbb{F}_{q^2}$.

A trinomial curve is the curve whose defining equation has precisely three monomials. The zeta function of the trinomial curve of type $k_1x^m+k_2y^n+1=0$ over a finite field $\mathbb{F}_{q}$ has been determined by Andr$\acute{e}$ Weil \cite{Weil} in 1949.

Some types of trinomial curves as maximal curves have been intensively studied. Gilles Lachaud \cite[Propostion 7]{Lachaud1} proved that the Fermat curve of equation $k_1x^m+k_2y^m+1=0,2\le m,gcd(q,m)=1$ defined over $\mathbb{F}_{q}$ is maximal over $\mathbb{F}_{q^2}$ if $m|(q+1)$. Angela Aguglia etc. \cite{Aguglia} proved that the Hurwitz curves of the form $x^ny+x+y^n=0$ with $gcd(q,n^2-n+1)=1$ is maximal over $\mathbb{F}_{q^2}$ if and only if $(n^2-n+1)|(q+1)$. They also proved that the curve of equation $x^ny^l+x^l+y^n=0$, where $n\ge l\ge 2$ and $gcd(q,n^2-nl+l^2)=1$, is maximal over $\mathbb{F}_{q^2}$ if $(n^2-nl+l^2)|(q+1)$. Arnaldo Garcia and Saeed Tafazolian \cite{Garcia2} proved that the Fermat curve given by $x^m+y^m=1$ with $2\le m,gcd(q,m)=1$ is maximal over $\mathbb{F}_{q^2}$ if and only if $m|(q+1)$. In \cite{Tafazolian2}, Saeed Tafazolian and Fernando Torres proved that the curve given by $x^n+y^m=1$ with $2\le m,2\le n,gcd(q,mn)=1$ is maximal over $\mathbb{F}_{q^2}$ if and only if both integers $m$ and $n$ divide $q + 1$. They \cite{Tafazolian3} subsequently proved that the curve given by $y^n=x^m+x$, where $2\le m,2\le n,gcd(q,(m-1)n)=1,n|m$, is maximal over $\mathbb{F}_{q^2}$ if and only if $(n(m-1))|(q+1)$.

In this paper, we determine the zeta functions of trinomial curves in terms of Gauss and Jacobi sums, and we obtain an explicit formula of the genus of a trinomial curve over any finite field, then we study the conditions for an trinomial curve to be a maximal curve over the finite field $\mathbb{F}_{q^2}$.

The main result of this paper is the following theorem. We only state the case 5 of the classification of trinomial curves. The other cases is similar to case 5. We first give some notations. Let $p$ be a prime number, $m=p^km'\in \mathbb{Z}, (p,m')=1,k\geq0$, we will denote $m'$ by $m_p$ and denote the $p$-adic valuation of $\mathbb{Q}$ by $v_p(\cdot)$. Let $\xi=(\xi_1,\xi_2)$ be a pair of rational numbers, we denote by $\mu(\xi)$ or $\mu(\xi_1,\xi_2)$ be the smallest positive integer such that $(q^{\mu(\xi)}-1)\xi_i\equiv0\pmod{1}$ for $i=1,2$.

\begin{theorem}
Let $C$ be the nonsingular model over $\mathbb{F}_{q}$ of the geometrically irreducible curve given by
\[
x^{m_1}y^{n_1}+k_1x^m+k_2y^n=0,
\]
where $k_1,k_2\in\mathbb{F}_{q}^*; m_1+n_1>m, m_1+n_1>n$, and $n_1 \geq m_1$, if $m_1=n_1$ then $n\geq m$. Let $d=gcd(m,n,m_1,n_1)$. Let $p=char(\mathbb{F}_{q})$. Suppose $p\nmid d$. Let $i(C)=\frac{m_1n+mn_1-mn-d_1-d_2-d_3}{2}+1$, where $d_1=gcd(m_1,n_1-n), d_2=gcd(n_1,m_1-m)$ and $d_3=gcd(m,n)$. Let $g(C)$ be the genus of $C$ over $\mathbb{F}_{q}$. Let $\xi=(\xi_1,\xi_2)$ be a pair of rational numbers such that
\[
\left\{\begin{array}{rcrl}
m_1\xi_1&+&(m_1-m)\xi_2&\equiv0\pmod{1}\\
(n_1-n)\xi_1&+&n_1\xi_2&\equiv0\pmod{1}
\end{array}\right.
\]
and
\[
\xi_i\not\equiv0\pmod{1},v_p(\xi_i)\geq0,\text{ for i=1,2, },
(\xi_1+\xi_2)\not\equiv0\pmod{1}.\textup{\quad$(1.1)$}
\]
Then
\begin{enumerate}
\item[(1)] If $i(C)=0$, then $g(C)=0$. If $i(C)>0$, then
\[
g(C)=\frac{(m_1n+mn_1-mn)_p-(d_1)_p-(d_2)_p-(d_3)_p}{2}+1.
\]

\item[(2)] The numerator of the zeta function of the curve $C$ over $\mathbb{F}_{q}$ is
\[
\begin{array}{l}
P_{C}(U)\\
=\prod\limits_{\xi}(1+\frac{1}{q^{\mu(\xi)}}\chi_{\xi_1}(k_2^{-1})\chi_{\xi_2}(k_1^{-1})
g(\psi,\chi_{\xi_1})g(\psi,\chi_{\xi_2})g(\psi,\chi_{-\xi_1-\xi_2})U^{\mu(\xi)}).
\end{array}
\]

\item[(3)] Suppose further that $((m_1n+m(n_1-n))/d)_p|(q^l+1)$ for some $l$. Then $\mu(\xi)$ is even and $\mu(\xi)=2(l,\mu(\xi))$. Let $\mu(\xi)=2\nu(\xi)$, then the numerator of the zeta function of the curve $C$ over $\mathbb{F}_{q}$ is
\[
P_{C}(U)=\prod_\xi(1+q^{\nu(\xi)}U^{\mu(\xi)}),
\]

the product in (2) and (3) both being taking over all pairs $\xi=(\xi_1,\xi_2)$ satisfying $(1.1)$ but taking only one representative for each set of pairs $(q^\rho\xi_1,q^\rho\xi_2)$ with $0\leq \rho < \mu(\xi)$.

\item[(4)] Suppose $((m_1n+m(n_1-n))/d)_p|(q+1)$, i.e. $l=1$ in (3). Then $C$ is maximal over $\mathbb{F}_{q^2}$. Conversely, if $C$ is maximal over $\mathbb{F}_{q^2}$ and $g(C)>0$, then $((m_1n+m(n_1-n))/d)_p|(q^2-1)$.

\end{enumerate}

\end{theorem}

Note that recently, Saeed Tafazolian and Fernando Torres \cite{Tafazolian4} also studied this Hurwitz type curves using Weierstrass semigroups and Serre's covering. They proved that \cite[Proposition 2.3]{Tafazolian4} the curve defined by $x^{m_1}y^{n_1}+x^m+y^n=0$ with $m,m_1,n,n_1\in \mathbb{N},m_1n+m(n_1-n)\ge 1, gcd(q, m_1n+m(n_1-n))=1$ is maximal over $\mathbb{F}_{q^2}$ if $(m_1n+m(n_1-n))|(q+1)$. They also showed that \cite[Theorem 2.9]{Tafazolian4}, if $m,m_1\ge 0, n\ge 2,m\equiv 1\pmod{n},\delta=m_1n+m(1-n)\ge 2,\Delta=\lfloor\frac{\delta}{n}\rfloor\ge 1,gcd(\Delta+1,n-1)=1$, then the curve defined by $x^{m_1}y+x^m+y^n=0$ is maximal over $\mathbb{F}_{q^2}$ if and only if $(m_1n+m(1-n))|(q+1)$.

The rest of this paper is organized as follows. In Section \ref{section2} trinomial curves are classified into 5 cases and we recall known results about the irreducibility and genus of a trinomial curve. At the end of this section we using Serre's covering \cite[Proposition 6]{Lachaud1} to give a simple criticism to find the finite field $\mathbb{F}_{q^2}$ such that a trinomial curve defined over $\mathbb{F}_{q}$ is maximal over $\mathbb{F}_{q^2}$. Section \ref{section3} introduces some preliminary details on Gauss and Jacobi sums and the number of solutions in $\mathbb{F}_q$ of the system of equations $x^{m_1}=a_1,x^{m_2}=a_2,\cdots,x^{m_r}=a_r$ which we will use to calculate the zeta function of a trinomial curve in Section \ref{section4}. We then use the zeta function of a trinomial curve to study when a trinomial curve is a maximal curve in Section \ref{section4}.

\section{Trinomial curves} \label{section2}
In this section we mainly recall some known results about the irreducibility and genus of a trinomial curve.

\begin{definition}
Let $C$ be an affine curve over a field K. Suppose the reduced equation $F$ of $C$ has exactly three monomials, we will call $C$ a trinomial curve.
\end{definition}

\begin{prop}\label{prop:irreducible_Beelen}(\cite[Proposition 2.10]{Beelen})
Let $K$ be a field, and $F(x,y)=\alpha x^a+\beta x^by^c+\gamma y^d\in K[x,y]$, where $a,b,c,d$ are nonnegative integers, and $\alpha\beta\gamma\ne 0$. Assume that $(a,0),(b,c)$ and $(0,d)$ are three distinct points. Then $F(x,y)$ is absolute irreducible over $K$ if and only if $ac+bd\ne ad$ and the characteristic of $K$ doesn't divide $gcd(a,b,c,d)$.
\end{prop}

\begin{prop} \label{prop:classify}
Let C be an geometrically irreducible trinomial curve over a field K, and let $p=char(K)$. Then the equation of $C$ is one of the following 5 cases with respect to permutation of variables:
\begin{enumerate}
%  \item  $k_1x^n+k_2y^n+1=0$, $p\nmid n$,
  \item  $k_1x^m+k_2y^n+1=0$, $n\ge m$, $p\nmid gcd(m,n)$,
  \item  $k_1x^m+k_2y^{n_1}+y^{n_2}=0$, $n_2>m, n_2>n_1$, $p\nmid gcd(m,n_1,n_2)$,
  \item  $k_1x^{m_1}y^{n_1}+k_2y^{n}+1=0$, $n>m_1+n_1$, $p\nmid gcd(m_1,n_1,n)$,
  \item  $k_1x^{m_1}y^{n_1}+k_2x^{m_2}y^{n_2}+1=0$, $m_1+n_1\geq m_2+n_2, \frac{n_1}{m_1} \geq \frac{n_2}{m_2}$, $(m_1,n_1)\ne (m_2,n_2)$, $p\nmid gcd(m_1,m_2,n_1,n_2)$,
  \item  $x^{m_1}y^{n_1}+k_1x^m+k_2y^n=0$, $m_1+n_1>m, m_1+n_1>n, n_1 \geq m_1$, and if $m_1=n_1$ then $n\geq m$. $p\nmid gcd(m_1,m,n_1,n)$,
\end{enumerate}
where $m,m_1,n,n_1$ are all positive integers, and $k_1,k_2\in K, k_1k_2\ne 0$.
\end{prop}
\proof
Suppose the homogeneous equation of the curve C is
$
k_1x^{a_1}y^{b_1}z^{c_1}+k_2x^{a_2}y^{b_2}z^{c_2}+k_3x^{a_3}y^{b_3}z^{c_3}=0.
$ We can divide it into two cases with respect to permutation of variables:
\begin{enumerate}
   \item $a_1=\min\{a_1,a_2,a_3\}$,$b_1=\min\{b_1,b_2,b_3\}$ and $c_2=\min\{c_1,c_2,c_3\}$. Then
   \[
   \begin{array}{l}
       k_1x^{a_1}y^{b_1}z^{c_1}+k_2x^{a_2}y^{b_2}z^{c_2}+k_3x^{a_3}y^{b_3}z^{c_3}\\
       =x^{a_1}y^{b_1}z^{c_2}(k_1z^{c_1-c_2}+k_2x^{a_2-a_1}y^{b_2-b_1}+k_3x^{a_3-a_1}y^{b_3-b_1}z^{c_3-c_2}).
   \end{array}
   \]
   \item $a_1=\min\{a_1,a_2,a_3\}$,$b_2=\min\{b_1,b_2,b_3\}$ and $c_3=\min\{c_1,c_2,c_3\}$. Then
   \[
   \begin{array}{l}
       k_1x^{a_1}y^{b_1}z^{c_1}+k_2x^{a_2}y^{b_2}z^{c_2}+k_3x^{a_3}y^{b_3}z^{c_3}\\
       =x^{a_1}y^{b_2}z^{c_3}(k_1y^{b_1-b_2}z^{c_1-c_3}+k_2x^{a_2-a_1}z^{c_2-c_3}+k_3x^{a_3-a1}y^{b_3-b_2}).
   \end{array}
   \]
 \end{enumerate}
Thus the homogeneous equation of the irreducible curve $C$ is of the form
   $
   k_1z^{a_1'}+k_2x^{a_2'}y^{b_2'}+k_3x^{a_3'}y^{b_3'}z^{c_3'}=0
   $ or
   $
   k_1x^{a_1'}y^{b_1'}+k_2y^{b_2'}z^{c_2'}+k_3x^{a_3'}z^{c_3'}=0.
   $
The former case can be divided into the following 6 reduced forms with respect to permutation of variables:
 \begin{enumerate}
   \item $a_2'=b_3'=c_3'=0$, $k_1z^{a_1'}+k_2y^{b_2'}+k_3x^{a_3'}=0,$
   \item $a_2'=b_3'=0$, $k_1z^{a_1'}+k_2y^{b_2'}+k_3x^{a_3'}z^{c_3'}=0,$
   \item $a_3'=0$, interchange $y$ and $z$, $k_1y^{a_1'}+k_2x^{a_2'}z^{b_2'}+k_3y^{c_3'}z^{b_3'}=0,$
   \item $a_2'=0$, $k_1z^{a_1'}+k_2y^{b_2'}+k_3x^{a_3'}y^{b_3'}z^{c_3'}=0,$
   \item $c_3'=0$, $(a_2',b_2')\ne(a_3',b_3')$, $k_1z^{a_1'}+k_2x^{a_2'}y^{b_2'}+k_3x^{a_3'}y^{b_3'}=0,$
   \item $k_1z^{a_1'}+k_2x^{a_2'}y^{b_2'}+k_3x^{a_3'}y^{b_3'}z^{c_3'}=0.$
 \end{enumerate}
And from the latter case we can get only one additional reduced form $k_1x^{a_1'}y^{b_1'}+k_2y^{b_2'}z^{c_2'}+k_3x^{a_3'}z^{c_3'}=0$.
Therefore, dehomogenization with respect to $z$ gives the following reduced equations of $C$ with respect to permutation of variables:
 \begin{enumerate}
%  \item  $k_1x^n+k_2y^n+1=0$,
  \item  $k_1x^m+k_2y^n+1=0$, $n\ge m$,
  \item  $k_1x^m+k_2y^{n_1}+y^{n_2}=0$, $n_2>m, n_2>n_1$,
  \item  $k_1x^{m_1}y^{n_1}+k_2y^{n}+1=0$, $n>m_1+n_1$,
  \item  $k_1x^{m_1}y^{n_1}+k_2x^{m_2}y^{n_2}+1=0$, $m_1+n_1\geq m_2+n_2, \frac{n_1}{m_1} \geq \frac{n_2}{m_2}$, $(m_1,n_1)\ne (m_2,n_2)$,
  \item  $x^{m_1}y^{n_1}+k_1x^m+k_2y^n=0$, $m_1+n_1>m, m_1+n_1>n, n_1 \geq m_1$, and if $m_1=n_1$ then $n\geq m$.
\end{enumerate}
The other statements follows by Proposition \ref{prop:irreducible_Beelen}.
\qed

\begin{definition}
A lattice polygon is a polygon with all vertices at points of the lattice of integers. A simple polygon is a closed polygonal chain of line segments in the plane which do not have points in common other than the common vertices of pairs of consecutive segments.
\end{definition}

\begin{theorem}[Pick](\cite[Theorem 3.1]{Olds})
Let $P$ be a simple lattice polygon. Then the area $A$ of polygon $P$ is given by the formula
\[
A= i+\frac{b}{2}-1.
\]
Where i is the number of lattice points inside P and b is the number of lattice points on the boundary of P including vertices.
\end{theorem}

\begin{corollary} \label{corollary:i(P)}
 Let $P$ be a simple lattice polygon with $n$ vertices\linebreak $P_1(a_1,b_1),P_2(a_2,b_2),\cdots,P_n((a_n,b_n))$ such that $\overrightarrow{P_1P_2\cdots P_nP_1}$ is the positively oriented boundary of $P$, where $a_i,b_i$ are integers. Write $P_{n+1}(a_{n+1},b_{n+1})=P_1(a_1,b_1)$, then the number $i(P)$ of lattice points inside $P$ is
 \[
 1+\sum_{i=1}^n\frac{a_ib_{i+1}-a_{i+1}b_i-d_{i,i+1}}{2},
 \]
 where $d_{i,i+1}=Gcd(a_i-a_{i+1},b_i-b_{i+1})$.
 Particularly, the number of lattice points inside a lattice triangle with vertices $P_1(a_1,b_1),P_2(a_2,b_2)$ and $P_3(a_3,b_3)$ is
 \[
 1+\frac{(a_1b_2-a_2b_1-d_1)+(a_2b_3-a_3b_2-d_2)+(a_3b_1-a_1b_3-d_3)}{2}
 \], where $d_1=Gcd(a_1-a_2,b_1-b_2),d_2=Gcd(a_2-a_3,b_2-b_3)$ and $d_3=Gcd(a_1-a_3,b_1-b_3)$.
\end{corollary}

\proof
The area of $P$ is $\sum\limits_{i=1}^n\frac{a_ib_{i+1}-a_{i+1}b_i}{2}$ by Green' theorem. The number of lattice points lying on the line segment $\overline{P_iP_{i+1}}$ including vertices $P_i,P_{i+1}$ is $d_{i,i+1}+1$. Therefore from Pick's theorem the number $i(P)$ of lattice points inside $P$ is
 \[
 1+\sum_{i=1}^n\frac{a_ib_{i+1}-a_{i+1}b_i-d_{i,i+1}}{2}.
 \]
\qed

\begin{definition}\cite{Beelen}
Let $K$ be a field. Let $F(x,y)=\sum_{i \in \mathcal{I}}\alpha_ix^{i_1}y^{i_2}$ be a polynomial in $K[x,y]$, with $i_1, i_2 \geq 0$. Denote by $\Gamma(F)$ the convex hull of the points $P_i=(i_1,i_2)$ in $\mathbb{R}_{\geq 0}^2$. The set $\Gamma(F)$ is called the Newton polygon of $F$. Let $C$ be a trinomial curve. We denote the Newton polygon of the equation $F$ of $C$ by $\Gamma(C)$ and denote the number of integral points in the interior of the Newton polygon $\Gamma(C)$ by $i(C)$.
\end{definition}

\begin{prop}\label{prop:i(C)}
Let $C$ be an absolute irreducible affine trinomial curve over a field K, and let $p=char(K)$. Then the value of $i(C)$ is listed as follows.
\begin{enumerate}
%  \item  $k_1x^n+k_2y^n+1=0$, $p\nmid n$,
  \item  $k_1x^m+k_2y^n+1=0$, $n\ge m$, $p\nmid gcd(m,n)$. $i(C)=\frac{(m-1)(n-1)-(d-1)}{2}$, where $d=gcd(m,n)$.
  \item  $k_1x^m+k_2y^{n_1}+y^{n_2}=0$, $n_2>m, n_2>n_1$, $p\nmid gcd(m,n_1,n_2)$. $i(C)=\frac{(m-1)(n_2-n_1)-(d_1+d_2)}{2}+1$, where $d_1=gcd(m,n_1),d_2=gcd(m,n_2)$.
  \item  $k_1x^{m_1}y^{n_1}+k_2y^{n}+1=0$, $n>m_1+n_1$, $p\nmid gcd(m_1,n_1,n)$. $i(C)=\frac{(m_1-1)n-d_1-d_2}{2}+1$, where $d_1=gcd(m_1,n_1), d_2=gcd(m_1,n-n_1)$.
  \item  $k_1x^{m_1}y^{n_1}+k_2x^{m_2}y^{n_2}+1=0$, $m_1+n_1\geq m_2+n_2, \frac{n_1}{m_1} \geq \frac{n_2}{m_2}$, $(m_1,n_1)\ne (m_2,n_2)$, $p\nmid gcd(m_1,m_2,n_1,n_2)$. $i(C)=\frac{m_2n_1-m_1n_2-d_1-d_2-d_3}{2}+1$, where $d_1=gcd(m_1,n_1),d_2=gcd(m_2,n_2)$ and $d_3=gcd(n_1-n_2,m_1-m_2)$.
  \item  $x^{m_1}y^{n_1}+k_1x^m+k_2y^n=0$, $m_1+n_1>m, m_1+n_1>n, n_1 \geq m_1$, and if $m_1=n_1$ then $n\geq m$. $p\nmid gcd(m_1,m,n_1,n)$. \linebreak $i(C)=\frac{m_1n+mn_1-mn-d_1-d_2-d_3}{2}+1$. where $d_1=gcd(m_1,n_1-n), d_2=gcd(n_1,m_1-m)$ and $d_3=gcd(m,n)$.
\end{enumerate}
\end{prop}
\proof
The value of $i(C)$ is obtained from Corollary \ref{corollary:i(P)}.
\qed

\begin{theorem}[Baker](\cite[Theorem 4.2]{Beelen})\label{thm:Baker}
Let $F(x,y)=0$ define an irreducible curve $\mathcal{X}$ over an algebraically closed field. Let $i(F)$ denote the number of integral points in the interior of the Newton polygon $\Gamma(F)$. Then the genus $g$ of the nonsingular model of $\mathcal{X}$ satisfies $g \leq i(F)$.  Equality holds if $F$ is nondegenerate with respect to its Newton polygon and the singular points of the homogeneous curve with equation $F^*(x,y,z)=0$ are among $[0:0:1],[0:1:0]$ and $[1:0:0]$.
\end{theorem}

\begin{prop}\cite[Corollary 4.3]{Beelen}
Let $\mathbb{F}$ be an algebraically closed field. Let $a,b,c$ and $d$ be nonnegative integers. Suppose a curve is given by the equation $\alpha x^a+\beta x^by^c+\gamma y^d=0$, where $\alpha,\beta,\gamma\in\mathbb{F}^*$, and $ac+bd\ne ad$, and $p$, the characteristic of $\mathbb{F}$ does not divide all of $a,b,c$ and $d$. Then the genus $g$ of the nonsingular model of this curve satisfies
\[
g\le 1+\frac{1}{2}\{|ac+bd-ad|-gcd(a-b,c)-gcd(b,c-d)-gcd(a,d)\}.
\]
If $p$ does not divide $gcd(a-b,c),gcd(b,c-d),gcd(a,d)$ and $ac+bd-ad$, then equality holds.
\end{prop}

\begin{prop}
Let $C$ be the nonsingular model over $\mathbb{F}_{q}$ of a geometrically irreducible curve with a homogeneous equation $k_1x^{a_{11}}y^{a_{12}}z^{a_{13}}+k_2x^{a_{21}}y^{a_{22}}z^{a_{23}}+x^{a_{31}}y^{a_{32}}z^{a_{33}}=0$, where $k_1,k_2\in\mathbb{F}_{q},k_1k_2\ne 0$. Let
\[
\mathbf{A} =
\left( \begin{array}{ccc}
a_{11} & a_{12} & a_{13}\\
a_{21} & a_{22} & a_{23}\\
a_{31} & a_{32} & a_{33}
\end{array} \right).
\]
Suppose the matrix $\mathbf{A}$ is nonsingular and n is the least positive integer such that two columns of $n\mathbf{A}^{-1}$ are all integers, then $C$ is maximal over finite fields $\mathbb{F}_{q^2}$ with $n|(q+1)$.
\end{prop}
\proof
Suppose
\[
n\mathbf{A}^{-1}=\left( \begin{array}{ccc}
b_{11} & b_{12} & b_{13}\\
b_{21} & b_{22} & b_{23}\\
b_{31} & b_{32} & b_{33}
\end{array} \right),
\]
and the first two columns are integers. let
$
\varphi(u,v)=(u^{b_{11}}v^{b_{12}},u^{b_{21}}v^{b_{22}}).
$
We can see that the Fermat curve of the form $C': k_1u^n+k_2v^n+1=0$ is a covering of the curve given by $k_1x^{a_{11}}y^{a_{12}}+k_2x^{a_{21}}y^{a_{22}}+x^{a_{31}}y^{a_{32}}=0$ by the morphism $\varphi$.
Since the curve $C'$ is maximal over finite fields $\mathbb{F}_{q^2}$ with $n|(q+1)$, we have $C$ is also maximal over the finite fields by \cite[Proposition 6]{Lachaud1}.
\qed

\begin{example}
The curve $xy^5+x^2y^3z+z^6=0$ over $\mathbb{F}_{q}$ with $7\nmid q$, which belongs to the case (4). Its genus is 3.  It's maximal over $\mathbb{F}_{q^2}$ with $7|(q+1)$.
\[
\mathbf{A}=\left( \begin{array}{ccc}
2 & 3 & 1\\
1 & 5 & 0\\
0 & 0 & 6
\end{array} \right),\quad
\mathbf{A}^{-1}=\left( \begin{array}{ccc}
\frac{5}{7} & -\frac{3}{7} & -\frac{5}{42}\\
-\frac{1}{7} & \frac{2}{7} & \frac{1}{42}\\
0 & 0 & \frac{1}{6}
\end{array} \right),
\]
\end{example}

\section{Exponential sum and the number of points} \label{section3}
In this section we recall some properties of Gauss and Jacobi sums and develop some preliminary results which will be used in Section \ref{section4} to calculate the zeta function of trinomial curves.

\begin{definition}(\cite[Chapter 4]{Katz})\label{def:gauss_sum}
Let $\chi:\mathbb{F}_q^*\rightarrow \mathbb{Q}(\zeta_{q-1})^*$ be a (possibly trivial) multiplicative character of $\mathbb{F}_q^*$ and $\psi:(\mathbb{F}_q, +)\rightarrow \mathbb{Q}(\zeta_p)^*$ a non-trivial additive character of $\mathbb{F}_q$, and let $a\in\mathbb{F}_q$. The Gauss sum is defined by
\[
g_a(\psi,\chi)=\sum_{u\in \mathbb{F}_q^*}\psi(au)\chi(u).
\]
For a fixed choice of non-trivial $\psi$ and any function $f$ on $\mathbb{F}_q^*$, say with values in an overfield $E$ of $\mathbb{Q}(\zeta_{q-1})$, we define its multiplicative Fourier transform $\hat{f}$ to be the $E$-valued function on characters given by
\[
\hat{f}(\chi)=\sum_{u\in \mathbb{F}_q^*}f(u)\chi(u).
\]
The Fourier inversion formula, $f(u)=\frac{1}{q-1}\sum_\chi\bar{\chi}(u)\hat{f}(\chi)$ allows us to recover $f$ from $\hat{f}$, where $u\in\mathbb{F}_q^*$. Given two functions $f,g$ on $\mathbb{F}_q^*$, let $u\in\mathbb{F}_q^*$, then their convolution $f*g$ is the function on $\mathbb{F}_q^*$ defined by
\[
(f*g)(u)=\sum_{xy=u}f(x)g(y).
\]
The Fourier transform of the convolution is given by the product of the transforms:
\[
(\hat{f*g})(\chi)=\hat{f}(\chi)\hat{g}(\chi).
\]
\end{definition}

\begin{remark}
Gauss sums are usually defined as (e.\ g.\ \cite{Berndt},\cite{Ireland},\cite{Weil})
\[
g_a'(\psi,\chi)=\sum_{u\in \mathbb{F}_q}\psi(au)\chi(u).
\]
These two definitions of Gauss sums are related as follows,
\[
g_a(\psi,\chi)=g_a'(\psi,\chi)-\chi(0)=\left\{ \begin{array}{ll}
g_a'(\psi,\chi) & \textrm{if $\chi$ is nontrivial}\\
g_a'(\psi,\chi)-1 & \textrm{if $\chi$ is trivial}.
\end{array}\right.
\]
\end{remark}

\begin{definition}\label{def:jacobi_sum}%(\cite{Berndt})
Let $\chi_1,\chi_2,\cdots,\chi_k$ be multiplicative characters on a finite field $\mathbb{F}_q^*$. A Jacobi sum is defined by the formula
\[
j(\chi_1,\chi_2,\cdots,\chi_k)=\sum_{t_1+t_2+\cdots+t_k=1}\chi_1(t_1)\chi_2(t_2)\cdots\chi_k(t_k).
\]
where the summation is taken over all $k$-tuples $(t_1,t_2,\cdots,t_k)$ of elements of $\mathbb{F}_q^*$ with $t_1+t_2+\cdots+t_k=1$.
$j_0(\chi_1,\chi_2,\cdots,\chi_k)$ is the character sum defined by
\[
j_0(\chi_1,\chi_2,\cdots,\chi_k)=\sum_{t_1+t_2+\cdots+t_k=0}\chi_1(t_1)\chi_2(t_2)\cdots\chi_k(t_k),
\]
where the summation is taken over all $k$-tuples $(t_1,t_2,\cdots,t_k)$ of elements of $\mathbb{F}_q^*$ with $t_1+t_2+\cdots+t_k=0$.
\end{definition}

\begin{remark}
Jacobi sum is usually defined as (e.\ g.\ \cite{Berndt},\cite{Ireland})
\[
j'(\chi_1,\chi_2,\cdots,\chi_k)=\sum_{t_1+t_2+\cdots+t_k=1 \atop t_i\in\mathbb{F}_q, i=1,2,\cdots,k}\chi_1(t_1)\chi_2(t_2)\cdots\chi_k(t_k)
\]
and $j_0'(\chi_1,\chi_2,\cdots,\chi_k)$ is defined as
\[
j_0'(\chi_1,\chi_2,\cdots,\chi_k)=\sum_{t_1+t_2+\cdots+t_k=0 \atop t_i\in\mathbb{F}_q, i=1,2,\cdots,k}\chi_1(t_1)\chi_2(t_2)\cdots\chi_k(t_k)
\]
It is easy to see if $\chi_1,\chi_2,\cdots,\chi_k$ are all nontrivial, then  $j(\chi_1,\chi_2,\cdots,\chi_k)=j'(\chi_1,\chi_2,\cdots,\chi_k)$ and $j_0(\chi_1,\chi_2,\cdots,\chi_k)=j_0'(\chi_1,\chi_2,\cdots,\chi_k)$.
\end{remark}

\begin{theorem} Let $\chi,\chi_1,\cdots,\chi_k$ be multiplicative characters on a finite field $\mathbb{F}_q^*$, $\psi$ be a fixed non-trivial additive character of $\mathbb{F}_q$ and $a\in \mathbb{F}_q$. Then
\begin{enumerate}
  \item
      \[
      j_0(\chi_1,\cdots,\chi_k)=\left\{\begin{array}{l}
      \frac{(q-1)^k+(-1)^k(q-1)}{q}\hspace{2em} \textrm{if $\chi_1,\cdots,\chi_k$ are all trivial.}\vspace{1em}\\
      (-1)^sj_0(\chi_{s+1},\chi_{s+2},\cdots,\chi_k)\hspace{1em} \textrm{if $k\geq 2$, and $s$ }\\
      \hspace{6em}\textrm{ characters are trivial, $1\leq s<k$,}\\
       \hspace{6em}\textrm{say $\chi_1,\chi_2,\cdots,\chi_s$ are all trivial}\vspace{1em}\\
      -(q-1)j(\chi_1,\cdots,\chi_k)\hspace{1em} \textrm{if $\chi_1,\cdots,\chi_k$ are all}\\
       \hspace{6.5em}\textrm{nontrivial and $\chi_1\cdots\chi_k$ is trivial.}\vspace{1em}\\
      0 \hspace{9em} \hspace{1em} \textrm{if $\chi_1\cdots\chi_k$ is nontrivial.}
      \end{array}\right.
      \]

  \item
      \[
      j(\chi_1,\cdots,\chi_k)=\left\{\begin{array}{l}
      \frac{(q-1)^k+(-1)^{k-1}}{q}\hspace{4em} \textrm{if $\chi_1,\cdots,\chi_k$ are all trivial.}\vspace{1em}\\
      (-1)^sj(\chi_{s+1},\chi_{s+2},\cdots,\chi_k)\hspace{1em} \textrm{if $k\geq 2$, and $s$ }\\
      \hspace{8em}\textrm{ characters are trivial, $1\leq s<k$,}\\
       \hspace{8em}\textrm{say $\chi_1,\chi_2,\cdots,\chi_s$ are all trivial}
      \end{array}\right.
      \]

  \item If $\chi_1,\cdots,\chi_k$ are nontrivial and $\chi_1\cdots\chi_k$ trivial, then
      \[
      \begin{array}{l}
      j_0(\chi_1,\cdots,\chi_k)=\chi_k(-1)(q-1)j(\chi_1,\cdots,\chi_{k-1}),\\
      j(\chi_1,\cdots,\chi_k)=-\chi_k(-1)j(\chi_1,\cdots,\chi_{k-1}).
      \end{array}
      \]
  \item If $\chi_1,\cdots,\chi_k$ are nontrivial, then
      \[
      j(\chi_1,\cdots,\chi_k)=\left\{\begin{array}{ll}
      \frac{g(\psi,\chi_1)\cdots g(\psi,\chi_k)}{g(\psi,\chi_1\cdots\chi_k)}& \textrm{if $\chi_1\cdots\chi_k$ is nontrivial.}\\
      \\
      -\frac{g(\psi,\chi_1)\cdots g(\psi,\chi_k)}{q}& \textrm{if $\chi_1\cdots\chi_k$ is trivial.}
      \end{array}\right.
      \]
  \item \[
      g_a(\psi,\chi)=\left\{\begin{array}{ll}
      \chi(a^{-1})g(\psi,\chi)&\textrm{ if $\chi$ is nontrivial and $a\neq 0$,}\\
      0&\textrm{ if $\chi$ is nontrivial and $a=0$,}\\
      -1&\textrm{ if $\chi$ is trivial and $a\neq 0$,}\\
      q-1&\textrm{ if $\chi$ is trivial and $a=0$,}
      \end{array}\right.
      \]
  \item \[
      \sum_{t\in\mathbb{F}_q^*}\chi(t)=\left\{\begin{array}{ll}
      0&\textrm{ if $\chi$ is nontrivial,}\\
      q-1&\textrm{ if $\chi$ is trivial.}
      \end{array}\right.
      \]
  \item If $a\in\mathbb{F}_q^*$. Then
    \[
      g_a(\psi,\chi)g_a(\psi,\bar{\chi})=\left\{\begin{array}{ll}
      \chi(-1)q&\textrm{ if $\chi$ is nontrivial,}\\
      1&\textrm{ if $\chi$ is trivial.}
      \end{array}\right.
    \]
\end{enumerate}
\end{theorem}

\proof
We only need to prove the first two case for the first two properties. For the details of the other properties see \cite{Berndt}.
\begin{enumerate}
\item Let $a_k:=j_0(\epsilon,\epsilon,\cdots,\epsilon)=\#\{(t_1,t_2,\cdots,t_k):t_1+t_2+\cdots+t_k=0,t_i\in \mathbb{F}_q^*\}$, where there are k trivial characters $\epsilon$'s in $j_0$. If $t_1,t_2,\cdots,t_{k-1}$ are chosen arbitrarily in $\mathbb{F}_q^*$, then $t_k$ is uniquely determined by the condition $t_1+t_2+\cdots+t_k=0$, but in order to $t_k\ne0$, we need $t_1+t_2+\cdots+t_{k-1}\ne0$. Thus for $k>1$,
    \[
    \begin{array}{rl}
    a_k=&(q-1)^{k-1}-\#\{(t_1,\cdots,t_{k-1}):t_1+t_2+\cdots+t_{k-1}=0,t_i\in \mathbb{F}_q^*\}\\
    =&(q-1)^{k-1}-a_{k-1}.
    \end{array}
    \]
    Note that $a_1=0,a_2=q-1$. Hence $a_k=\frac{(q-1)^k+(-1)^k(q-1)}{q}$.

    For the second case, we first show that if $\chi_k=\epsilon,\chi_{k-1}\ne\epsilon$, then
    \[
    j_0(\chi_1,\cdots,\chi_{k-1},\chi_k)=-j_0(\chi_1,\cdots,\chi_{k-1}).
    \]
    In fact,
    \[
     \begin{array}{l}
    j_0(\chi_1,\cdots,\chi_{k-1},\epsilon)\\
    =\sum\limits_{t_1+\cdots+t_{k-1}+t_k=0\atop t_i\in\mathbb{F}_q^*, i=1,\cdots,k-1,k}\chi_1(t_1)\cdots\chi_{k-1}(t_{k-1})\\
    =\sum\limits_{t_1+\cdots+t_{k-2}=0,t_{k-1}=-t_k\atop t_i\in\mathbb{F}_q^*, i=1,\cdots,k-1,k}\chi_1(t_1)\cdots\chi_{k-1}(t_{k-1})\\
    \qquad+\sum\limits_{a\in\mathbb{F}_q^*}\sum\limits_{t_1+\cdots+t_{k-2}=a,t_{k-1}\ne-a\atop t_i\in\mathbb{F}_q^*, i=1,\cdots,k-1,k}\chi_1(t_1)\cdots\chi_{k-1}(t_{k-1})\\
    =j_0(\chi_1,\cdots,\chi_{k-2})\cdot\sum\limits_{t_{k-1}\in\mathbb{F}_q^*}\chi_{k-1}(t_{k-1})+\\
    \quad\sum\limits_{a\in\mathbb{F}_q^*}\sum\limits_{t_1+\cdots+t_{k-2}=a\atop t_i\in\mathbb{F}_q^*}\chi_1(t_1)\cdots\chi_{k-2}(t_{k-2})\cdot\sum\limits_{t_{k-1}\ne-a\atop t_{k-1}\in\mathbb{F}_q^*}\chi_{k-1}(t_{k-1})\\
    =\sum\limits_{a\in\mathbb{F}_q^*}\sum\limits_{t_1+\cdots+t_{k-2}=a\atop t_i\in\mathbb{F}_q^*}\chi_1(t_1)\cdots\chi_{k-2}(t_{k-2})(-\chi_{k-1}(-a))\\
    =-j_0(\chi_1,\cdots,\chi_{k-1})
    \end{array}
    \]
    Note that the value of $j_0$ isn't affected by the order of $\chi_i$'s. Hence if there are s characters are trivial and at least one character is nontrivial, then we can iterate the above result for s times and obtain the required result.
\item Let $b_k=j(\epsilon,\epsilon,\cdots,\epsilon)=\#\{(t_1,\cdots,t_{k-1},t_k):t_1+\cdots+t_{k-1}+t_k=1,t_i\in \mathbb{F}_q^*\}$. The condition for $t_i$ can be divided into the following two cases. One is $t_1+\cdots+t_{k-1}=0,t_k=1,t_i\in \mathbb{F}_q^*,1\le i\le k$, the other is $t_1+\cdots+t_{k-1}=a,t_k=1-a,a\ne 0,1,t_i\in \mathbb{F}_q^*,1\le i\le k$, which implies that $\frac{t_1}{a}+\cdots+\frac{t_{k-1}}{a}=1$. Thus for $k>1$, $b_k=a_{k-1}+(q-2)b_{k-1}$. Note that $b_1=1,b_2=q-2$. Hence $b_k=\frac{(q-1)^k+(-1)^{k-1}}{q}$.

    The second case for $j$ is similar to that of $j_0$, so we omit the proof.
\end{enumerate}

\qed

\begin{corollary}\label{corollary:jj0}
Let $\chi_1,\chi_2,\cdots,\chi_k$ be multiplicative characters on $\mathbb{F}_q^*$. Let $\chi_{k+1}$\linebreak$=\overline{\chi_1\chi_2\cdots\chi_k}$. Then
\[
j(\chi_1,\chi_2,\cdots,\chi_k)=\frac{1}{q-1}\chi_{k+1}(-1)j_0(\chi_1,\chi_2,\cdots,\chi_k,\chi_{k+1}).
\]
\end{corollary}
\proof
First, suppose $\chi_1,\cdots,\chi_k$ are all trivial, then $\chi_{k+1}$ is trivial, and \linebreak $j(\chi_1,\cdots,\chi_k)=\frac{(q-1)^k+(-1)^{k-1}}{q}$, hence
\[
j_0(\chi_1,\cdots,\chi_k,\chi_{k+1})=\frac{(q-1)^{k+1}+(-1)^{k+1}(q-1)}{q}=(q-1)j(\chi_1,\cdots,\chi_k).
\]

Second, suppose $\chi_1,\cdots,\chi_k$ are all nontrivial. If $\chi_{k+1}$ is nontrivial, then this corollary is true by the above theorem. If $\chi_{k+1}$ is trivial, i.e. $\chi_{k+1}=\overline{\chi_1\chi_2\cdots\chi_k}=\epsilon$, then $\overline{\chi_1\chi_2\cdots\chi_{k-1}}=\chi_k\ne\epsilon$. Thus from the former case we have
\[
\begin{array}{rl}
j(\chi_1,\cdots,\chi_k)&=-\chi_k(-1)j(\chi_1,\cdots,\chi_{k-1})\\
&=-\chi_k(-1)\frac{1}{q-1}\chi_k(-1)j_0(\chi_1,\cdots,\chi_{k-1},\chi_k)\\
&=-\frac{1}{q-1}j_0(\chi_1,\cdots,\chi_k)\\
&=\frac{1}{q-1}j_0(\chi_1,\cdots,\chi_k,\epsilon).
\end{array}
\]

Finally, suppose there are s characters nontrivial and k-s characters trivial, say $\chi_{s+1}=\cdots=\chi_k=\epsilon$, then $\chi_{k+1}=\overline{\chi_1\chi_2\cdots\chi_s}$. Thus
\[
\begin{array}{rl}
j(\chi_1,\cdots,\chi_s,\epsilon,\cdots,\epsilon)&=(-1)^{k-s}j(\chi_1,\cdots,\chi_s)\\
&=(-1)^{k-s}\frac{1}{q-1}\chi_{k+1}(-1)j_0(\chi_1,\cdots,\chi_s,\chi_{k+1})\\
&=\frac{1}{q-1}\chi_{k+1}(-1)j_0(\chi_1,\cdots,\chi_s,\epsilon,\cdots,\epsilon,\chi_{k+1}).
\end{array}
\]

\qed

From the above theorem, let $\chi_\beta,\chi_{\beta_1},\chi_{\beta_2},\chi_{\beta_3}$ be multiplicative characters on a finite field $\mathbb{F}_q^*$ and $\psi$ be a fixed non-trivial additive character of $\mathbb{F}_q$. Then we have
\[
j_0(\bar{\chi}_{\beta_1},\bar{\chi}_{\beta_2},\bar{\chi}_{\beta_3})=\left\{\begin{array}{l}
      (q-1)(q-2) \hspace{5em} \textrm{if $\chi_{\beta_1},\chi_{\beta_2},\chi_{\beta_3}$ are all trivial,}\vspace{1em}\\
      \frac{(q-1)}{q}g(\psi,\bar{\chi}_{\beta_1})g(\psi,\bar{\chi}_{\beta_2})g(\psi,\bar{\chi}_{\beta_3})
       \hspace{1em} \textrm{if $\chi_{\beta_1},\chi_{\beta_2},\chi_{\beta_3}$ are}\\
       \hspace{7em}\textrm{not all trivial and $\chi_{\beta_1+\beta_2+\beta_3}$ is trivial,}\vspace{1em}\\
      0 \hspace{12em}\textrm{otherwise.}
      \end{array}\right.
\]
\[
g_a(\psi,\chi_\beta)=\left\{\begin{array}{ll}
      \chi_\beta(a^{-1})g(\psi,\chi_\beta)& \textrm{if $a\neq 0$,}\\
      q-1& \textrm{if $\chi_\beta$ is trivial and $a=0$,}\\
      0& \textrm{if $\chi_\beta$ is nontrivial and $a=0$,}
      \end{array}\right.
\]

\begin{lemma}
Let $\chi$ be multiplicative characters on a finite field $\mathbb{F}_q^*$ and $\psi$ be a fixed non-trivial additive character of $\mathbb{F}_q$, and let $b \in \mathbb{F}_q^*$. Then
\[
\sum_{a\in\mathbb{F}_q}g_{ab}(\psi,\chi)=0.
\]
\end{lemma}

\proof
\[
\begin{array}{rl}
\sum\limits_{a\in\mathbb{F}_q}g_{ab}(\psi,\chi)&=g_0(\psi,\chi)+\sum\limits_{a\in\mathbb{F}_q^*}g_{ab}(\psi,\chi)\\
&=(q-1)\chi(0)+g(\psi,\chi)\sum\limits_{a\in\mathbb{F}_q^*}\chi(a^{-1}b^{-1})\\
&=(q-1)\chi(0)+g(\psi,\chi)(q-1)\chi(0)\\
&=0
\end{array}
\]

\qed

\begin{lemma}
Let $\chi_1,\chi_2$ be multiplicative characters on a finite field $\mathbb{F}_q^*$ and $\psi$ be a fixed non-trivial additive character of $\mathbb{F}_q$, and let $b_1,b_2 \in \mathbb{F}_q^*$. Then
\[
\sum_{a\in\mathbb{F}_q}g_{ab_1}(\psi,\chi_1)g_{ab_2}(\psi,\chi_2)=q(q-1)\chi_1(-b_1^{-1}b_2)\chi_1\chi_2(0).
\]
\end{lemma}

\proof
\[
\begin{array}{l}
\sum\limits_{a\in\mathbb{F}_q}g_{ab_1}(\psi,\chi_1)g_{ab_2}(\psi,\chi_2)\\
=g_0(\psi,\chi_1)g_0(\psi,\chi_2)
+\sum\limits_{a\in\mathbb{F}_q^*}g_{ab_1}(\psi,\chi_1)g_{ab_2}(\psi,\chi_2)\\
=(q-1)^2\chi_1(0)\chi_2(0)+g(\psi,\chi_1)g(\psi,\chi_2)\chi_1(b_1^{-1})\chi_2(b_2^{-1})
\sum\limits_{a\in\mathbb{F}_q^*}\chi_1\chi_2(a^{-1})\\
=(q-1)^2\chi_1(0)\chi_2(0)+g(\psi,\chi_1)g(\psi,\chi_2)\chi_1(b_1^{-1})\chi_2(b_2^{-1})(q-1)\chi_1\chi_2(0)\\

=\left\{\begin{array}{ll}
q(q-1)\chi_1(-b_1^{-1}b_2)&\text{if $\chi_1\ne\epsilon$ and $\chi_1\chi_2=\epsilon$}\\
q(q-1)&\text{if $\chi_1=\chi_2=\epsilon$}\\
0&\text{if $\chi_1\chi_2\ne\epsilon$}\\
\end{array}\right.\\

=q(q-1)\chi_1(-b_1^{-1}b_2)\chi_1\chi_2(0)\\
\end{array}
\]

\qed

From now on we usually use (m,n) to denote the g.\ c.\ d.\ of m and n.

\begin{definition}
Let $\mathbb{F}_q^*=<\zeta>$, $\chi$ be a multiplicative character of $\mathbb{F}_q^*$ with $\chi(\zeta)=e^{\frac{2\pi ik}{q-1}}=e^{2\pi i\alpha}$, (we will denote such a $\chi$ by $\chi_{\alpha,\zeta}$ or $\chi_{\alpha}$ below), $m$ be a positive integer, $(m, q-1)=d$, and $ms+(q-1)t=d$. If $d|k$, then $\chi^{\frac{1}{m}}$ is defined by
\[
\chi^{\frac{1}{m}}(\zeta)=e^{\frac{2\pi ik}{(q-1)}\frac{s}{d}}.
\]
In fact, $\chi^{\frac{1}{m}}$ is one of the multiplicative character $\chi_0$ such that $\chi_0^m=\chi$ in the multiplicative character group over $\mathbb{F}_q^*$.
\end{definition}

\begin{lemma}\label{lemma:uni_solu}
Let $\chi_\delta$ be multiplicative characters over $\mathbb{F}_q^*$. Suppose $(m, q-1)=d, ms+(q-1)t=d$. Then
\[
\begin{array}{l}
(\chi_\alpha\chi_\beta)^{\frac{1}{m}}\chi_\gamma=\chi_\delta,  \frac{q-1}{d}(\alpha+\beta)\equiv0\pmod{1}, d\gamma\equiv0\pmod{1}\\
\hspace{14em}\Leftrightarrow (\beta,\gamma)\equiv(m\delta-\alpha,\frac{(q-1)t}{d}\delta)\pmod{1}.
\end{array}
\]
\end{lemma}
\proof
$(\chi_\alpha\chi_\beta)^{\frac{1}{m}}\chi_\gamma=\chi_\delta$ is equivalent to $\frac{s}{d}(\alpha+\beta)+\gamma\equiv\delta\pmod{1}$, and the latter is equivalent to
$\frac{s(q-1)}{d}(\alpha+\beta)+(q-1)\gamma\equiv(q-1)\delta\pmod{q-1}$. From
$d\gamma\equiv0\pmod{1}$, i.e. $(q-1)\gamma\equiv0\pmod{\frac{q-1}{d}}$, we have
$\frac{s(q-1)}{d}(\alpha+\beta)\equiv(q-1)\delta\pmod{\frac{q-1}{d}}$, which implies that
$\frac{(q-1)}{d}(\alpha+\beta)\equiv\frac{(q-1)m}{d}\delta\pmod{\frac{q-1}{d}}$, hence $\beta\equiv m\delta-\alpha\pmod{1}$.
Substitute $\beta\equiv m\delta-\alpha\pmod{1}$ to $\frac{s}{d}(\alpha+\beta)+\gamma\equiv\delta\pmod{1}$, we obtain $\gamma\equiv(1-\frac{ms}{d})\delta\equiv\frac{(q-1)t}{d}\delta\pmod{1}$.

\qed

\begin{prop} \label{prop:exponential_Gauss}
Suppose $\mathbb{F}_q^*=<\zeta>$, $\chi$ is a multiplicative character of $\mathbb{F}_q^*$ with $\chi(\zeta)=e^{\frac{2\pi ik}{q-1}}=e^{2\pi i\alpha}$, and $\psi$ is a fixed non-trivial additive character of $\mathbb{F}_q$. Let $a\in\mathbb{F}_q$, and let $m$ be a positive integer and $(m, q-1)=d$. Then
\[
\sum_{u\in\mathbb{F}_q^*}\chi(u)\psi(au^m)=\left\{ \begin{array}{ll}
0 & \textrm{if $d\nmid (q-1)\alpha$}\\
\sum\limits_\beta g_a(\psi,\chi^{\frac{1}{m}}\chi_\beta) & \textrm{if $d|(q-1)\alpha$,}
\end{array} \right.
\]
where $d\beta\equiv0\pmod{1}$.
\end{prop}

\proof
First suppose $(m,q-1)=1$. Let $u^m=w\in \mathbb{F}_q^*$. Then $u=w^{\frac{1}{m}}$. Hence
\[
\sum_{u\in\mathbb{F}_q^*}\chi(u)\psi(au^m)=\sum_{w\in\mathbb{F}_q^*}\chi(w^{\frac{1}{m}})\psi(aw)
=\sum_{w\in\mathbb{F}_q^*}\chi^{\frac{1}{m}}(w)\psi(aw)
=g_a(\psi,\chi^{\frac{1}{m}})
\]

Second, suppose $m|(q-1)$ and $m|k$. Let $u^m=w\in \mathbb{F}_q^*$, and $\theta$ be a primitive mth root of unity in $\mathbb{F}_q^*$. Then $w^{\frac{1}{m}},w^{\frac{1}{m}}\theta,w^{\frac{1}{m}}\theta^2,\cdots,w^{\frac{1}{m}}\theta^{m-1}$ are all the solutions of
$x^m=w$ in $\mathbb{F}_q^*$. Let $N_m(a)$ denote the number of solutions in $\mathbb{F}_q^*$ of the equation $x^m=a$ with $a\in\mathbb{F}_q^*$. Then
\[
\sum_{u\in\mathbb{F}_q^*}\chi(u)\psi(au^m)=\sum_{w\in\mathbb{F}_q^*\atop i=0,1,\cdots,m-1}\chi(w^{\frac{1}{m}}\theta^i)\psi(aw)\frac{N_m(w)}{m}
\]
Note that $\chi(w^{\frac{1}{m}}\theta^i)=\chi^{\frac{1}{m}}((w^{\frac{1}{m}}\theta^i)^m)=\chi^{\frac{1}{m}}(w)$.
It follows that
\[
\begin{array}{ll}
\sum\limits_{u\in\mathbb{F}_q^*}\chi(u)\psi(au^m)=\sum\limits_{w\in\mathbb{F}_q^*}\chi^{\frac{1}{m}}(w)\psi(aw)N_m(w)
&=\sum\limits_\beta\sum_{w\in\mathbb{F}_q^*}\chi^{\frac{1}{m}}\chi_\beta(w)\psi(aw)\\
&=\sum\limits_\beta g_a(\psi,\chi^{\frac{1}{m}}\chi_\beta)
\end{array}
\]
where $m\beta\equiv0\pmod{1}$.

Third, suppose $m|(q-1)$ and $m\nmid k$. Then
\[
\begin{array}{ll}
\sum\limits_{u\in\mathbb{F}_q^*}\chi(u)\psi(au^m)&=\sum\limits_{w\in\mathbb{F}_q^*\atop i=0,1,\cdots,m-1}\chi(w^{\frac{1}{m}}\theta^i)\psi(aw)\frac{N_m(w)}{m}\\
&=(\sum\limits_{i=0}^{m-1}\chi(\theta)^i)\sum\limits_{w\in\mathbb{F}_q^*}\chi(w^{\frac{1}{m}})\psi(aw)\frac{N_m(w)}{m}\\
&=0
\end{array}
\]
because $\chi(\theta)\ne 1$.

Finally we prove the general case by induction on $m$. Suppose $d>1$. Let $m_1=\frac{m}{d},d_1=(m_1,q-1)$, and $\theta$ be a primitive dth root of unity in $\mathbb{F}_q^*$. Then $m_1<m$ and $d_1|d$. Thus
\[
\begin{array}{ll}
\sum\limits_{u\in\mathbb{F}_q^*}\chi(u)\psi(au^m)&=\sum\limits_{u\in\mathbb{F}_q^*}\chi(u)\psi(a(u^d)^{m_1})\\
&=\sum\limits_{w\in\mathbb{F}_q^*\atop i=0,1,\cdots,d-1}\chi(w^{\frac{1}{d}}\theta^i)\psi(aw^{m_1})\frac{N_d(w)}{d}\\
&=\left\{\begin{array}{ll}
0 & \textrm{if $d\nmid (q-1)\alpha$}\\
\\
\sum\limits_{w\in\mathbb{F}_q^*}\chi(w^{\frac{1}{d}})\psi(aw^{m_1})N_d(w)
& \textrm{if $d|(q-1)\alpha$}\\
=\sum\limits_\beta\sum_{w\in\mathbb{F}_q^*}\chi^{\frac{1}{d}}\chi_\beta(w)\psi(aw^{m_1})\\
\end{array}\right.
\end{array}
\]
where $d\beta\equiv0\pmod{1}$.
Now let $\alpha_1=\frac{\alpha}{d}+\beta$. Then by induction
\[
\sum_{w\in\mathbb{F}_q^*}\chi^{\frac{1}{d}}\chi_\beta(w)\psi(aw^{m_1})
=\left\{ \begin{array}{ll}
0 & \textrm{if $d_1\nmid (q-1)\alpha_1$}\\
\sum\limits_{\beta_1}\sum\limits_{v\in\mathbb{F}_q^*}(\chi^{\frac{1}{d}}\chi_\beta)^{\frac{1}{m_1}}\chi_{\beta_1}(v)\psi(av) & \textrm{if $d_1|(q-1)\alpha_1$}
\end{array} \right.,
\]
where $d_1\beta_1\equiv0\pmod{1}$.

Now suppose $(m,q-1)=d, ms+(q-1)t=d, (m_1,q-1)=d_1, m_1s_1+(q-1)t_1=d_1$. We will show that if $d|(q-1)\alpha$, then the sets $\{(\chi^{\frac{1}{d}}\chi_\beta)^{\frac{1}{m_1}}\chi_{\beta_1}:d\beta\equiv0\pmod{1},d_1|(q-1)(\frac{\alpha}{d}+\beta),
d_1\beta_1\equiv0\pmod{1}\}$ and $\{\chi^{\frac{1}{m}}\chi_\gamma:d\gamma\equiv0\pmod{1}\}$ are the same, that is, the sets $A=\{((\frac{\alpha}{d}+\beta)\frac{s_1}{d_1}+\beta_1)\pmod{1}:d\beta\equiv0\pmod{1},d_1|(q-1)(\frac{\alpha}{d}+\beta),
d_1\beta_1\equiv0\pmod{1}\}$ and $B=\{(\frac{\alpha}{d}s+\gamma)\pmod{1}:d\gamma\equiv0\pmod{1}\}$ are the same.

From $(\frac{q-1}{d},d_1)=1$ as $d_1|m_1$, we can see there exist $\frac{d}{d_1}$ solutions for $\beta \pmod{1}$ satisfying equations $d\beta\equiv0\pmod{1},d_1|(q-1)(\frac{\alpha}{d}+\beta)$. And if the pairs $(\beta,\beta_1)$ are distinct, then the characters $(\chi^{\frac{1}{d}}\chi_\beta)^{\frac{1}{m_1}}\chi_{\beta_1}$ are also distinct by Lemma \ref{lemma:uni_solu}. Thus $\#A=\#B=d$. Moreover, since $\frac{m_1}{d_1}s_1\equiv1\pmod{\frac{q-1}{d_1}}$ and $ \frac{q-1}{d}|\frac{q-1}{d_1}$, $\frac{m_1}{d_1}s_1\equiv1\pmod{\frac{q-1}{d}}$. From $m_1\equiv d_1\frac{m_1}{d_1}\pmod{\frac{q-1}{d}}$, we have $(\frac{m_1}{d_1})^{-1}\equiv d_1m_1^{-1}\pmod{\frac{q-1}{d}}$, which implies that $s_1\equiv d_1s\pmod{\frac{q-1}{d}}$. Note that if $d_1\beta_1\equiv0\pmod{1}$, then $(q-1)\beta_1\equiv0\pmod{\frac{q-1}{d_1}}$, hence $(q-1)\beta_1\equiv0\pmod{\frac{q-1}{d}}$. Thus given any pair $(\beta,\beta_1)$ such that $d\beta\equiv0\pmod{1},d_1|(q-1)(\frac{\alpha}{d}+\beta),
d_1\beta_1\equiv0\pmod{1}$, we have $(q-1)((\frac{\alpha}{d}+\beta)\frac{s_1}{d_1}+\beta_1-\frac{\alpha}{d}s)\equiv0\pmod{\frac{q-1}{d}}$, then $d((\frac{\alpha}{d}+\beta)\frac{s_1}{d_1}+\beta_1-\frac{\alpha}{d}s)\equiv0\pmod{1}$, that is, $\gamma\equiv((\frac{\alpha}{d}+\beta)\frac{s_1}{d_1}+\beta_1-\frac{\alpha}{d}s)\pmod{1}$ is the only one solution for $\gamma \pmod{1}$ satisfying equations $(\frac{\alpha}{d}+\beta)\frac{s_1}{d_1}+\beta_1\equiv\frac{\alpha}{d}s+\gamma\pmod{1},d\gamma\equiv0\pmod{1}$. Therefore $A=B$, i.e.
\[
\begin{array}{ll}
\sum\limits_{u\in\mathbb{F}_q^*}\chi(u)\psi(au^m)=\left\{\begin{array}{ll}
0 & \textrm{if $d\nmid (q-1)\alpha$}\\
\\
\sum\limits_\beta\sum\limits_{\beta_1}\sum\limits_{v\in\mathbb{F}_q^*}
(\chi^{\frac{1}{d}}\chi_\beta)^{\frac{1}{m_1}}\chi_{\beta_1}(v)\psi(av)
& \textrm{if $d|(q-1)\alpha$}\\
=\sum\limits_\gamma\sum\limits_{v\in\mathbb{F}_q^*}\chi^{\frac{1}{m}}\chi_\gamma(v)\psi(av)\\
\end{array}\right.
\end{array}
\]
where $d\gamma\equiv0\pmod{1}$.
\qed

\begin{remark}
The book\cite[1.1.4]{Berndt} gives a formula $\sum_{u\in\mathbb{F}_q^*}\psi(au^k)=\sum_{j=1}^{k-1}g_a(\psi,\linebreak\chi^j)$, where $\chi$ is a character of order $k$ on $\mathbb{F}_q$ and $a\in\mathbb{F}_q^*$, which is only valid for $k|(q-1)$. However, the above proposition is already known to I.\ M.\ Gel'fand, etc. \cite[9.2]{Gelfand} as a corollary of their Proposition 8.5.
\end{remark}

\begin{corollary}
Let $\chi_\alpha, \chi_\beta$ be multiplicative character of $\mathbb{F}_q^*$ and $\psi$ be a fixed non-trivial additive character of $\mathbb{F}_q$. Let $a\in\mathbb{F}_q$, and let $m$ be a positive integer and $(m, q-1)=d$. Then
\[
\sum_{u\in\mathbb{F}_q^*}\chi_\alpha(u)\chi_\beta(u)\psi(au^m)=\left\{ \begin{array}{ll}
0 & \textrm{if $d\nmid (q-1)(\alpha+\beta)$}\\
\sum\limits_{\gamma}g_a(\psi,(\chi_\alpha\chi_\beta)^{\frac{1}{m}}\chi_\gamma) & \textrm{if $d|(q-1)(\alpha+\beta)$}
\end{array} \right.,
\]
where $d\gamma\equiv0 \pmod{1}$.
\end{corollary}

\qed

In other words, let $f_1(t)=\psi(at^m),f_2(t)=\chi_\alpha(t)\psi(at^m)$ and $(m, q-1)=d$. Then
\[
\hat{f}_1(\chi_\beta)=\sum_{u\in\mathbb{F}_q^*}\chi_\beta(u)\psi(au^m)=\left\{ \begin{array}{ll}
0 & \textrm{if $d\nmid (q-1)\beta$}\\
\sum\limits_{\gamma}g_a(\psi,\chi_\beta^{\frac{1}{m}}\chi_\gamma) & \textrm{if $d|(q-1)\beta$}
\end{array} \right.,
\]
where $d\gamma\equiv0\pmod{1}$.
\[
\begin{array}{ll}
\hat{f}_2(\chi_\beta)&=\sum_{u\in\mathbb{F}_q^*}\chi_\beta(u)\chi_\alpha(u)\psi(au^m)\\
&=\left\{ \begin{array}{ll}
0 & \textrm{if $d\nmid (q-1)(\beta+\alpha)$}\\
\sum\limits_{\gamma}g_a(\psi,(\chi_\beta\chi_\alpha)^{\frac{1}{m}}\chi_\gamma)& \textrm{if $d|(q-1)(\beta+\alpha)$,}
\end{array} \right.
\end{array}
\]
where $d\gamma\equiv0\pmod{1}$.

Hence the corresponding Fourier inversion formulae are
\[
\begin{array}{ll}
f_1(t)=\psi(at^m)=\frac{1}{q-1}\sum\limits_\beta\bar{\chi}_\beta(t)\hat{f}_1(\chi_\beta)
&=\frac{1}{q-1}\sum\limits_\beta\bar{\chi}_\beta(t)\sum\limits_\gamma g_a(\psi,\chi_\beta^{\frac{1}{m}}\chi_\gamma)\\
&=\frac{1}{q-1}\sum\limits_{\beta,\gamma}\bar{\chi}_\beta(t)g_a(\psi,\chi_\beta^{\frac{1}{m}}\chi_\gamma),
\end{array}
\]
where $\frac{q-1}{d}\beta\equiv0\pmod{1}, d\gamma\equiv0\pmod{1}$.

\[
\begin{array}{ll}
f_2(t)=\chi_\alpha(t)\psi(at^m)&=\frac{1}{q-1}\sum\limits_\beta\bar{\chi}_\beta(t)\hat{f}_2(\chi_\beta)\\
&=\frac{1}{q-1}\sum\limits_\beta\bar{\chi}_\beta(t)\sum\limits_\gamma g_a(\psi,(\chi_\alpha\chi_\beta)^{\frac{1}{m}}\chi_\gamma)\\
&=\frac{1}{q-1}\sum\limits_{\beta,\gamma}\bar{\chi}_\beta(t)g_a(\psi,(\chi_\alpha\chi_\beta)^{\frac{1}{m}}\chi_\gamma),
\end{array}
\]
where $\frac{q-1}{d}(\alpha+\beta)\equiv0\pmod{1}, d\gamma\equiv0\pmod{1}$.

Now let $f=\chi_{\alpha_1}(t)\psi(a_1t^{m_1}),g=\chi_{\alpha_2}(t)\psi(a_2t^{m_2}), d_1=(m_1, q-1), d_2=(m_2,q-1)$. Then
\[
(f*g)(t)=\sum_{uv=t}f(u)g(v)=\sum_{uv=t}\chi_{\alpha_1}(u)\psi(a_1u^{m_1})\chi_{\alpha_2}(v)\psi(a_2v^{m_2}).
\]
Since $(\hat{f*g})(\chi)=\hat{f}(\chi)\hat{g}(\chi)$, we have
\[
\begin{array}{ll}
(\hat{f*g})(\chi_\beta)
&=\sum\limits_{t\in\mathbb{F}_q^*}\chi_\beta(t)(\sum\limits_{uv=t}\chi_{\alpha_1}(u)\psi(a_1u^{m_1})\chi_{\alpha_2}(v)\psi(a_2v^{m_2}))\\
&=(\sum\limits_{t\in\mathbb{F}_q^*}\chi_\beta(t)\chi_{\alpha_1}(t)\psi(a_1t^{m_1}))(\sum\limits_{t\in\mathbb{F}_q^*}\chi_\beta(t)\chi_{\alpha_2}(t)\psi(a_2t^{m_2}))\\
&=\left\{ \begin{array}{l}
0 \qquad \textrm{if $d_1\nmid (q-1)(\alpha_1+\beta)$ or $d_2\nmid (q-1)(\alpha_2+\beta)$}\\
\\
\sum\limits_{\gamma_1,\gamma_2}g_{a_1}(\psi,(\chi_\beta\chi_{\alpha_1})^{\frac{1}{m_1}}\chi_{\gamma_1})
 g_{a_2}(\psi,(\chi_\beta\chi_{\alpha_2})^{\frac{1}{m_2}}\chi_{\gamma_2})
\end{array} \right.
\\
& \qquad\qquad\quad \textrm{if $d_1|(q-1)(\alpha_1+\beta)$ and $d_2|(q-1)(\alpha_2+\beta)$},
\end{array}
\]
where $d_1\gamma_1\equiv0\pmod{1}, d_2\gamma_2\equiv0\pmod{1}$.

The corresponding Fourier inversion formula is
\[
\begin{array}{ll}
(f*g)(t)&=\sum\limits_{uv=t}\chi_{\alpha_1}(u)\psi(a_1u^{m_1})\chi_{\alpha_2}(v)\psi(a_2v^{m_2})\\
&=\frac{1}{q-1}\sum\limits_\beta\bar{\chi}_\beta(t)(\hat{f*g})(\chi_\beta)\\
&=\frac{1}{q-1}\sum\limits_{\beta,\gamma_1,\gamma_2}\bar{\chi}_\beta(t)g_{a_1}(\psi,(\chi_\beta\chi_{\alpha_1})^{\frac{1}{m_1}}\chi_{\gamma_1})
 g_{a_2}(\psi,(\chi_\beta\chi_{\alpha_2})^{\frac{1}{m_2}}\chi_{\gamma_2}),
\end{array}
\]
where $\beta,\gamma_1,\gamma_2$ satisfy $\frac{q-1}{d_1}(\alpha_1+\beta)\equiv0\pmod{1}, d_1\gamma_1\equiv0\pmod{1}$ and $\frac{q-1}{d_2}(\alpha_2+\beta)\equiv0\pmod{1}, d_2\gamma_2\equiv0\pmod{1}$.

\begin{definition}
Let $\mathbb{F}_q$ be a finite field. We denote the number of solutions in $\mathbb{F}_q$ of the equation $x^m=a$ by $N_{m}(a)$ and denote the number of solutions in $\mathbb{F}_q$ of the system of equations $x^{m_1}=a_1,x^{m_2}=a_2,\cdots,x^{m_r}=a_r$ by \linebreak $N_{m_1,m_2,\cdots,m_r}(a_1,a_2,\cdots,a_r)$, i.e.
\[
N_{m_1,m_2,\cdots,m_r}(a_1,a_2,\cdots,a_r)=\sharp\{x\in\mathbb{F}_q|x^{m_1}=a_1,x^{m_2}=a_2,\cdots,x^{m_r}=a_r\}.
\]
\end{definition}

\begin{prop}
Let $m_1,m_2,\cdots,m_r$ be positive integers. Suppose\linebreak $(m_1,m_2,\cdots,m_r)=d,m_1s_1+m_2s_2+\cdots+m_rs_r=d$ and $m_i=m_i'd$. Then
\[
\begin{array}{l}
 N_{m_1,m_2,\cdots,m_r}(a_1,a_2,\cdots,a_r)\\
 =N_d(a_1^{s_1}a_2^{s_2}\cdots a_r^{s_r})\cdot\delta(a_1^{m_2's_2+m_3's_3+\cdots+m_r's_r},a_2^{m_1's_2}a_3^{m_1's_3}\cdots a_r^{m_1's_r})\\
 \quad\cdot \delta(a_2^{m_1's_1+m_3's_3+\cdots+m_r's_r},a_1^{m_2's_1}a_3^{m_2's_3}\cdots a_r^{m_2's_r})\\
 \quad\cdots\\
 \quad\cdot \delta(a_r^{m_1's_1+m_2's_2+\cdots+m_{r-1}'s_{r-1}},a_1^{m_r's_1}a_2^{m_r's_2}\cdots a_{r-1}^{m_r's_{r-1}}),
\end{array}
 \]
where $\delta(x, y)=1$ if $x=y$ and zero otherwise.
\end{prop}
\proof
If $x^{m_1}=a_1,x^{m_2}=a_2,\cdots,$ and $x^{m_r}=a_r$, then
\[
x^d=(x^{m_1})^{s_1}(x^{m_2})^{s_2}\cdots(x^{m_r})^{s_r}=a_1^{s_1}a_2^{s_2}\cdots a_r^{s_r}
\]
and
\[
\begin{array}{ll}
a_1^{m_2's_2+m_3's_3+\cdots+m_r's_r}&=a_1^{\frac{m_2s_2+m_3s_3+\cdots+m_rs_r}{d}}\\
&=(x^{m_1})^{\frac{m_2s_2+m_3s_3+\cdots+m_rs_r}{d}}\\
&=(x^{m_1'})^{m_2s_2+m_3s_3+\cdots+m_rs_r}\\
&=a_2^{m_1's_2}a_3^{m_1's_3}\cdots a_r^{m_1's_r}.
\end{array}
\]
Similarly for others. Conversely, if $x^d=a_1^{s_1}a_2^{s_2}\cdots a_r^{s_r}$, then the equality \linebreak $a_1^{m_2's_2+m_3's_3+\cdots+m_r's_r}=a_2^{m_1's_2}a_3^{m_1's_3}\cdots a_r^{m_1's_r}$ is equivalent to the equality\linebreak
$x^{m_1}=a_1$ since
\[
\begin{array}{ll}
x^{m_1}=(x^d)^{m_1'}=(a_1^{s_1}a_2^{s_2}\cdots a_r^{s_r})^{m_1'}
&=a_1^{m_1's_1}a_2^{m_1's_2}a_3^{m_1's_3}\cdots a_r^{m_1's_r}\\
&=a_1^{m_1's_1}a_1^{m_2's_2+m_3's_3+\cdots+m_r's_r}\\
&=a_1.
\end{array}
\]
Similarly for others.

\qed

\begin{corollary}
Let $m,n$ be positive integers. Suppose $(m,n)=d,ms+nt=d,m=m_1d,n=n_1d$. Then $N_{m,n}(a,b)=N_d(a^sb^t)\delta(a^{n_1},b^{m_1})$, where $\delta(x, y)=1$ if $x=y$ and zero otherwise.
\end{corollary}

\proof
From the above proposition, we have
\[
N_{m,n}(a,b)=N_d(a^sb^t)\delta(a^{n_1t},b^{m_1t})\delta(b^{m_1s},a^{n_1s}).
\]
If $a^{n_1t}=b^{m_1t}$ and $b^{m_1s}=a^{n_1s}$, then $a^{n_1}=(a^{n_1})^{(sm_1+tn_1)}=(a^{n_1s})^{m_1}(a^{n_1t})^{n_1}=(b^{m_1s})^{m_1}(b^{m_1t})^{n_1}=(b^{m_1})^{(sm_1+tn_1)}=b^{m_1}$, as $m_1s+n_1t=1$.
\qed

\begin{prop}(\cite[Proposition 10.3.3]{Ireland})
Let $\alpha,x,y\in \mathbb{F}_q$. Then
\[
\frac{1}{q}\sum_{\alpha \in \mathbb{F}_q}\psi(\alpha(x-y))=\delta(x,y),
\]
where $\delta(x, y)=1$ if $x=y$ and zero otherwise.
\end{prop}

\section{Zeta function}\label{section4}

This section we employ the results developed in the former section to study trinomial curves. We will first calculate the zeta function of a trinomial curve $C$ over a finite field $\textbf{F}_q$, then determine the finite fields $\mathbb{F}_{q'^2}$ such that $C$ is maximal over these finite fields.

We consider case 5 before other cases.

\subsection{Case 5}

(5)$x^{m_1}y^{n_1}+k_1x^m+k_2y^n=0$    and $m_1+n_1>m, m_1+n_1>n$, and $n_1 \geq m_1$, if $m_1=n_1$ then $n\geq m$.

Suppose $C$ is an irreducible affine curve defined over a finite field $\mathbb{F}_q$ with equation of the form (5), i.\ e.\ $x^{m_1}y^{n_1}+k_1x^m+k_2y^n=0$    and $m_1+n_1>m, m_1+n_1>n$, and $n_1 \geq m_1$, if $m_1=n_1$ then $n\geq m$.  We will first express the number of points of the curve $C$ over $\mathbb{F}_q$ as sum of product of Gauss sum and Jacobi sum, then to compute the zeta function of the curve.

Let $(m_1,m)=d_1,m_1=m_1'd_1,m=m'd_1,m_1s_1+mt_1=d_1,(n_1,n)=d_2,n_1=n_1'd_2,n=n'd_2,n_1s_2+nt_2=d_2,(d_1,q-1)=d_1',(d_2,q-1)=d_2', (m_1',q-1)=l_1,(n_1',q-1)=l_2,(m',q-1)=l_3$ and $(n',q-1)=l_4$. Let $N$ be the number of points in $\mathbb{F}_q$ of the curve $C$. Put $L(u,v)=u_0v_0+k_1u_1+k_2v_1$. Then
\[
\begin{array}{lll}

N&=\sum\limits_{L(u,v)=0}N_{m_1,m}(u_0,u_1)N_{n_1,n}(v_0,v_1)&\\

&=\sum\limits_{L(u,v)=0}N_{d_1}(u_0^{s_1}u_1^{t_1})\delta(u_0^{m'},u_1^{m_1'})
N_{d_2}(v_0^{s_2}v_1^{t_2})\delta(v_0^{n'},v_1^{n_1'})&

\end{array}
\]
Now we divide $u,v$ into 4 cases:
\begin{enumerate}
\item $u_0\neq 0, u_1\neq 0, v_0\neq 0, v_1\neq 0$;
\item $u_0= 0, u_1=0, v_0\neq 0, v_1\neq 0$; which is impossible since then $L(u,v)\neq 0$.
\item $u_0\neq 0, u_1\neq 0, v_0=0, v_1=0$; which is also impossible for the same reason.
\item $u_0=0, u_1=0, v_0=0, v_1=0$;
\end{enumerate}

Thus we have
\[
\begin{array}{ll}
&N\\
&=1+\sum\limits_{L(u,v)=0\atop u,v\in\mathbb{F}_q^*}(\sum\limits_{\alpha_1}\chi_{\alpha_1}(u_0^{s_1}u_1^{t_1}))
(\sum\limits_{\alpha_2}\chi_{\alpha_2}(v_0^{s_2}v_1^{t_2}))(\frac{1}{q}\sum\limits_{a_1\in\mathbb{F}_q}\psi(a_1(u_0^{m'}-u_1^{m_1'})))\\

&\hspace{20em}\cdot(\frac{1}{q}\sum\limits_{a_2\in\mathbb{F}_q}\psi(a_2(v_0^{n'}-v_1^{n_1'})))\\
&\hspace{5em}(d_i'\alpha_i\equiv0\pmod{1},0\leq\alpha_i<1, \textrm{for i=1,2})\\

&=1+\frac{1}{q^2}\sum\limits_{L(u,v)=0\atop u,v\in\mathbb{F}_q^*}\sum\limits_{\alpha,a}\chi_{\alpha_1}(u_0^{s_1}u_1^{t_1})
\chi_{\alpha_2}(v_0^{s_2}v_1^{t_2})\psi(a_1(u_0^{m'}-u_1^{m_1'}))\\
&\hspace{18em}\cdot\psi(a_2(v_0^{n'}-v_1^{n_1'}))\\
&\hspace{5em}(d_i'\alpha_i\equiv0\pmod{1},a_i\in\mathbb{F}_q,0\leq\alpha_i<1, \textrm{for i=1,2})\\

&=1+\frac{1}{q^2}\sum\limits_{\alpha,a}\sum\limits_{L(u,v)=0\atop u,v\in\mathbb{F}_q^*}\chi_{\alpha_1}(u_1^{t_1})\psi(-a_1u_1^{m_1'})
\chi_{\alpha_2}(v_1^{t_2})\psi(-a_2v_1^{n_1'})\\
&\hspace{8em}\cdot\chi_{\alpha_1}(u_0^{s_1})\psi(a_1u_0^{m'})
\chi_{\alpha_2}(v_0^{s_2})\psi(a_2v_0^{n'})\\
% &\hspace{5em}(d_i'\alpha_i\equiv0\pmod{1},a_i\in\mathbb{F}_q,0\leq\alpha_i<1, \textrm{for i=1,2})\\

&=1+\frac{1}{q^2}\sum\limits_{\alpha,a}\sum\limits_{w_1+u_1+v_1=0\atop u_1,v_1,w_1\in\mathbb{F}_q^*}
\chi_{t_1\alpha_1}(k_1^{-1})\chi_{t_1\alpha_1}(u_1)\psi(-a_1k_1^{-m_1'}u_1^{m_1'})
\chi_{t_2\alpha_2}(k_2^{-1})\\
&\hspace{1em}
\cdot\chi_{t_2\alpha_2}(v_1)\psi(-a_2k_2^{-n_1'}v_1^{n_1'})\cdot\sum\limits_{u_0v_0=w_1\atop w_1\in\mathbb{F}_q^*}\chi_{s_1\alpha_1}(u_0)\psi(a_1u_0^{m'})
\chi_{s_2\alpha_2}(v_0)\psi(a_2v_0^{n'})\\
% &\hspace{5em}(d_i'\alpha_i\equiv0\pmod{1},a_i\in\mathbb{F}_q,0\leq\alpha_i<1, \textrm{for i=1,2})\\

&=1+\frac{1}{q^2}\sum\limits_{\alpha,a}\chi_{t_1\alpha_1}(k_1^{-1})\chi_{t_2\alpha_2}(k_2^{-1})\sum\limits_{u_1+v_1+w_1=0\atop u_1,v_1,w_1\in\mathbb{F}_q^*}\chi_{t_1\alpha_1}(u_1)\psi(-a_1k_1^{-m_1'}u_1^{m_1'})\\
&\hspace{1em}\cdot\chi_{t_2\alpha_2}(v_1)\psi(-a_2k_2^{-n_1'}v_1^{n_1'})
\cdot\sum\limits_{u_0v_0=w_1\atop w_1\in\mathbb{F}_q^*}\chi_{s_1\alpha_1}(u_0)\psi(a_1u_0^{m'})
\chi_{s_2\alpha_2}(v_0)\psi(a_2v_0^{n'})\\
% &\hspace{5em}(d_i'\alpha_i\equiv0\pmod{1},a_i\in\mathbb{F}_q,0\leq\alpha_i<1, \textrm{for i=1,2})\\
\end{array}
\]

\[
\begin{array}{ll}
&=1+\frac{1}{q^2}\sum\limits_{\alpha,a}\chi_{t_1\alpha_1}(k_1^{-1})\chi_{t_2\alpha_2}(k_2^{-1})\\
&\hspace{1em}\cdot\sum\limits_{u_1+v_1+w_1=0\atop u_1,v_1,w_1\in\mathbb{F}_q^*}
(\frac{1}{q-1}\sum\limits_{\beta_1,\gamma_1}\bar{\chi}_{\beta_1}(u_1)g_{-a_1k_1^{-m_1'}}
(\psi,(\chi_{t_1\alpha_1}\chi_{\beta_1})^{\frac{1}{m_1'}}\chi_{\gamma_1}))\\
&\hspace{1em}\cdot(\frac{1}{q-1}\sum\limits_{\beta_2,\gamma_2}\bar{\chi}_{\beta_2}(v_1)g_{-a_2k_2^{-n_1'}}
(\psi,(\chi_{t_2\alpha_2}\chi_{\beta_2})^{\frac{1}{n_1'}}\chi_{\gamma_2}))\\
&\hspace{1em}\cdot(\frac{1}{q-1}\sum\limits_{\beta_3,\gamma_3,\gamma_4}\bar{\chi}_{\beta_3}(w_1)
g_{a_1}(\psi,(\chi_{s_1\alpha_1}\chi_{\beta_3})^{\frac{1}{m'}}\chi_{\gamma_3})
g_{a_2}(\psi,(\chi_{s_2\alpha_2}\chi_{\beta_3})^{\frac{1}{n'}}\chi_{\gamma_4}))\\
&\hspace{5em}(a_i\in\mathbb{F}_q,d_i'\alpha_i\equiv0\pmod{1},0\leq\alpha_i<1, \textrm{for i=1,2}\\
&\hspace{5em}\frac{q-1}{l_1}(t_1\alpha_1+\beta_1)\equiv0\pmod{1},l_1\gamma_1\equiv0\pmod{1},\\
&\hspace{5em}\frac{q-1}{l_2}(t_2\alpha_2+\beta_2)\equiv0\pmod{1},l_2\gamma_2\equiv0\pmod{1},\\
&\hspace{5em}\frac{q-1}{l_3}(s_1\alpha_1+\beta_3)\equiv0\pmod{1},l_3\gamma_3\equiv0\pmod{1},\\
&\hspace{5em}\frac{q-1}{l_4}(s_2\alpha_2+\beta_3)\equiv0\pmod{1},l_4\gamma_4\equiv0\pmod{1}).(*)\\
%&\hspace{5em}\textrm{From now on, we will denote $g_a(\psi,\chi)$ just by $g_a(\chi)$}.&\\
\end{array}
\]

\[
\begin{array}{ll}
&=1+\frac{1}{q^2(q-1)^3}\sum\limits_{\alpha,\beta,\gamma,a}\chi_{t_1\alpha_1}(k_1^{-1})\chi_{t_2\alpha_2}(k_2^{-1})
g_{-a_1k_1^{-m_1'}}(\psi,(\chi_{t_1\alpha_1}\chi_{\beta_1})^{\frac{1}{m_1'}}\chi_{\gamma_1})\\
&\hspace{3em}\cdot g_{-a_2k_2^{-n_1'}}(\psi,(\chi_{t_2\alpha_2}\chi_{\beta_2})^{\frac{1}{n_1'}}\chi_{\gamma_2})
\cdot g_{a_1}(\psi,(\chi_{s_1\alpha_1}\chi_{\beta_3})^{\frac{1}{m'}}\chi_{\gamma_3})\\
&\hspace{3em}\cdot g_{a_2}(\psi,(\chi_{s_2\alpha_2}\chi_{\beta_3})^{\frac{1}{n'}}\chi_{\gamma_4})
\cdot\sum\limits_{u_1+v_1+w_1=0\atop u_1,v_1,w_1\in\mathbb{F}_q^*}\bar{\chi}_{\beta_1}(u_1)\bar{\chi}_{\beta_2}(v_1)\bar{\chi}_{\beta_3}(w_1)\\

&=1+\frac{1}{q^2(q-1)^3}\sum\limits_{\alpha,\beta,\gamma,a}\chi_{t_1\alpha_1}(k_1^{-1})\chi_{t_2\alpha_2}(k_2^{-1})
g_{-a_1k_1^{-m_1'}}(\psi,(\chi_{t_1\alpha_1}\chi_{\beta_1})^{\frac{1}{m_1'}}\chi_{\gamma_1})\\
&\hspace{3em}\cdot g_{-a_2k_2^{-n_1'}}(\psi,(\chi_{t_2\alpha_2}\chi_{\beta_2})^{\frac{1}{n_1'}}\chi_{\gamma_2})
\cdot g_{a_1}(\psi,(\chi_{s_1\alpha_1}\chi_{\beta_3})^{\frac{1}{m'}}\chi_{\gamma_3})\\
&\hspace{3em}\cdot g_{a_2}(\psi,(\chi_{s_2\alpha_2}\chi_{\beta_3})^{\frac{1}{n'}}\chi_{\gamma_4})
\cdot j_0(\bar{\chi}_{\beta_1},\bar{\chi}_{\beta_2},\bar{\chi}_{\beta_3})\\

&=1+\frac{1}{q-1}\sum\limits_{\alpha,\beta,\gamma}
\chi_{t_1\alpha_1}(k_1^{-1})\chi_{t_2\alpha_2}(k_2^{-1})((\chi_{t_1\alpha_1}\chi_{\beta_1})^{\frac{1}{m_1'}}\chi_{\gamma_1})
(k_1^{m_1'})\\
&\hspace{3em}\cdot((\chi_{t_2\alpha_2}\chi_{\beta_2})^{\frac{1}{n_1'}}\chi_{\gamma_2})(k_2^{n_1'})
j_0(\bar{\chi}_{\beta_1},\bar{\chi}_{\beta_2},\bar{\chi}_{\beta_3})\\

&\textrm{\quad where apart from the conditions $(*)$, the sum satisfies two additional}\\
&\textrm{conditions: }\\
&\hspace{8em}(\chi_{t_1\alpha_1}\chi_{\beta_1})^{\frac{1}{m_1'}}\chi_{\gamma_1}(\chi_{s_1\alpha_1}\chi_{\beta_3})^{\frac{1}{m'}}\chi_{\gamma_3}=\epsilon,\\
&\hspace{8em}(\chi_{t_2\alpha_2}\chi_{\beta_2})^{\frac{1}{n_1'}}\chi_{\gamma_2}(\chi_{s_2\alpha_2}\chi_{\beta_3})^{\frac{1}{n'}}\chi_{\gamma_4}=\epsilon\\

\\

&=1+\frac{1}{q-1}\sum\limits_{\alpha,\beta,\gamma}
\chi_{\beta_1}(k_1)\chi_{\beta_2}(k_2)j_0(\bar{\chi}_{\beta_1},\bar{\chi}_{\beta_2},\bar{\chi}_{\beta_3})\\

\end{array}
\]

Now let $(\chi_{t_1\alpha_1}\chi_{\beta_1})^{\frac{1}{m_1'}}\chi_{\gamma_1}=\chi_{\delta_1},
(\chi_{t_2\alpha_2}\chi_{\beta_2})^{\frac{1}{n_1'}}\chi_{\gamma_2}=\chi_{\delta_2}$. Then from the lemma \ref{lemma:uni_solu} we have
\[
\begin{array}{l}
N=1+\frac{1}{q-1}\sum\limits_{\alpha,\delta}\chi_{t_1\alpha_1-m_1'\delta_1}(k_1^{-1})\chi_{t_2\alpha_2-n_1'\delta_2}(k_2^{-1})\\
\hspace{8em}\cdot j_0(\chi_{t_1\alpha_1-m_1'\delta_1},\chi_{t_2\alpha_2-n_1'\delta_2},\chi_{s_1\alpha_1+m'\delta_1})
\end{array}
\]
where $\alpha,\delta$ satisfy
\[
\begin{array}{c}
s_1\alpha_1+m'\delta_1\equiv s_2\alpha_2+n'\delta_2\pmod{1},\\
d_i'\alpha_i\equiv0\pmod{1},(q-1)\delta_i\equiv0\pmod{1},\\
0\leq\alpha_i<1,0\leq\delta_i<1, \textrm{for i=1,2}
\end{array}
\]
and additionally
\[
t_1\alpha_1-m_1'\delta_1+t_2\alpha_2-n_1'\delta_2+s_1\alpha_1+m'\delta_1\equiv 0\pmod{1}
\]
due to the property of $j_0$.

Clearly $\alpha_i$ also satisfies $(q-1)\alpha_i\equiv0\pmod{1}$, hence $d_i'\alpha_i\equiv0\pmod{1}\Leftrightarrow d_i\alpha_i\equiv0\pmod{1}$, so we need to consider the following system of linear \linebreak congruence
\[
\left\{\begin{array}{rrrrr}
m'\delta_1&-n'\delta_2&+s_1\alpha_1&-s_2\alpha_2&\equiv 0\pmod{1}\\
(m'-m_1')\delta_1&-n_1'\delta_2&+(t_1+s_1)\alpha_1&+t_2\alpha_2&\equiv 0\pmod{1}\\
&&d_1\alpha_1&&\equiv0\pmod{1}\\
&&&d_2\alpha_2&\equiv0\pmod{1}
\end{array}\right.
\]

Let
\[
A=\left(\begin{array}{ccccc}
m'     &-n'   &s_1    &-s_2\\
m'-m_1'&-n_1' &t_1+s_1&t_2 \\
0      &0     &d_1    &0   \\
0      &0     &0      &d_2
\end{array}\right),
X=\left(\begin{array}{c}
\delta_1\\
\delta_2\\
\alpha_1\\
\alpha_2
\end{array}\right).
\]
 Then we have $AX\equiv0\pmod{1}$.

Following the method of \cite{Butson} or \cite{Smith}, we can show that
\[
\begin{array}{l}
\left(\begin{array}{cccc}
1     &0    &0    &0 \\
-1    &1    &0    &0 \\
0     &-m_1 &1    &0 \\
n_1   &-n_1 &0    &1
\end{array}\right)
A
\left(\begin{array}{cccc}
t_1  &t_1   &-s_1  &t_1 \\
0    &t_2   &s_2   &t_2+s_2\\
m_1' &m_1'  &m'    &m_1'   \\
0    &n_1'  &-n'   &n_1'-n'
\end{array}\right)\\
=
\left(\begin{array}{cccc}
1     &0   &0   &0 \\
0     &1   &0   &0 \\
0     &0   &m   &m_1\\
0     &0   &-n  &n_1-n
\end{array}\right)
\end{array}
\]
We denote the above equality by $UAV=B$, then $AX\equiv0\pmod{1}$$\Leftrightarrow$ \linebreak $ UAVV^{-1}X\equiv0\pmod{1}$. Let $V^{-1}X=Y=(y_1,y_2,y_3,y_4)^t$. Then $BY\equiv0\pmod{1}$, and from $X=VY$ we have $t_1\alpha_1-m_1'\delta_1=y_3,t_2\alpha_2-n_1'\delta_2=-y_3-y_4$ and $s_1\alpha_1+m'\delta_1=y_1+y_2+y_4$. Let $y_3=\theta_2, -y_3-y_4=\theta_1$. Then
\[
N=1+\frac{1}{q-1}\sum_{\theta_1,\theta_2}\chi_{\theta_1}(k_2^{-1})\chi_{\theta_2}(k_1^{-1})
j_0(\chi_{\theta_1},\chi_{\theta_2},\chi_{-\theta_1-\theta_2}),
\]
where $\theta_1,\theta_2$ satisfy $(q-1)\theta_i\equiv0\pmod{1}$ for $i=1,2$ and
\[
\begin{array}{l}
\left(\begin{array}{cc}
m_1   & m_1-m \\
n_1-n & n_1
\end{array}\right)
\left(\begin{array}{c}
\theta_1\\
\theta_2
\end{array}\right)
\equiv0\pmod{1}.
\end{array}
\]

For $j_0(\chi_{\theta_1},\chi_{\theta_2},\chi_{-\theta_1-\theta_2})$, we divide it into 5 cases.
\begin{enumerate}
  \item
  $\theta_1\equiv\theta_2\equiv0 \pmod{1}$, then $j_0(\chi_{\theta_1},\chi_{\theta_2},\chi_{-\theta_1-\theta_2})=(q-1)(q-2)$.
  \item
  $\theta_1\equiv0\pmod{1},\theta_2\not\equiv0\pmod{1}$, then $j_0(\chi_{\theta_1},\chi_{\theta_2},\chi_{-\theta_1-\theta_2})=-(q-1)\chi_{\theta_2}(-1)$, and from
  $ m_1\theta_1+(m_1-m)\theta_2\equiv0\pmod{1},(n_1-n)\theta_1+n_1\theta_2\equiv0\pmod{1}$, we have $(m_1-m)\theta_2\equiv0\pmod{1},n_1\theta_2\equiv0\pmod{1}$ and $(q-1)\theta_2\equiv0\pmod{1}$, which are equivalent to $(m_1-m,n_1,q-1)\theta_2\equiv0\pmod{1}$ and $\theta_2\not\equiv0\pmod{1}$, where $(m_1-m,n_1,q-1)$ denotes the g.c.d. of $m_1-m,n_1,q-1$.
  \item
  $\theta_2\equiv0\pmod{1},\theta_1\not\equiv0\pmod{1}$, then $j_0(\chi_{\theta_1},\chi_{\theta_2},\chi_{-\theta_1-\theta_2})=-(q-1)\chi_{\theta_1}(-1)$, and from
  $m_1\theta_1+(m_1-m)\theta_2\equiv0\pmod{1},(n_1-n)\theta_1+n_1\theta_2\equiv0\pmod{1}$, we have $m_1\theta_1\equiv0\pmod{1},(n_1-n)\theta_1\equiv0\pmod{1}$ and $(q-1)\theta_1\equiv0\pmod{1}$, which are equivalent to $(m_1,n_1-n,q-1)\theta_1\equiv0\pmod{1}$ and $\theta_1\not\equiv0\pmod{1}$.
  \item
  $\theta_1+\theta_2\equiv0\pmod{1},\theta_1\not\equiv0\pmod{1}$, then $j_0(\chi_{\theta_1},\chi_{\theta_2},\chi_{-\theta_1-\theta_2})=-(q-1)\chi_{\theta_1}(-1)$, and from
  $ m_1(\theta_1+\theta_2)-m\theta_2\equiv0\pmod{1},(n_1-n)(\theta_1+\theta_2)+n\theta_2\equiv0\pmod{1}$ and $\theta_2\equiv-\theta_1\pmod{1}$, we have $m\theta_1\equiv0\pmod{1},n\theta_1\equiv0\pmod{1}$ and $(q-1)\theta_1\equiv0\pmod{1}$, which are equivalent to $(m,n,q-1)\theta_1\equiv0\pmod{1}$ and $\theta_1\not\equiv0\pmod{1}$.
  \item
  $\theta_1\not\equiv0\pmod{1},\theta_2\not\equiv0\pmod{1},\theta_1+\theta_2\not\equiv0\pmod{1}$, then\linebreak
  $j_0(\chi_{\theta_1},\chi_{\theta_2},\chi_{-\theta_1-\theta_2})=\frac{q-1}{q}g(\psi,\chi_{\theta_1})g(\psi,\chi_{\theta_2})g(\psi,\chi_{-\theta_1-\theta_2})$.
\end{enumerate}

So
\[
\begin{array}{ll}
N=&q-1-\sum\limits_{(m_1-m,n_1,q-1)\xi_1\equiv0\pmod{1}\atop\xi_1\not\equiv0\pmod{1}}\chi_{\xi_1}(-k_1^{-1})\\
&-\sum\limits_{(m_1,n_1-n,q-1)\xi_2\equiv0\pmod{1}\atop\xi_2\not\equiv0\pmod{1}}\chi_{\xi_2}(-k_2^{-1})\\
&-\sum\limits_{(m,n,q-1)\xi_3\equiv0\pmod{1}\atop\xi_3\not\equiv0\pmod{1}}\chi_{\xi_3}(-k_1k_2^{-1})\\
&+\frac{1}{q}\sum\limits_{\zeta_1,\zeta_2}\chi_{\zeta_1}(k_2^{-1})\chi_{\zeta_2}(k_1^{-1})
g(\psi,\chi_{\zeta_1})g(\psi,\chi_{\zeta_2})g(\psi,\chi_{-\zeta_1-\zeta_2}),
\end{array}
\]
where in the last sum $\zeta_1,\zeta_2$ satisfy $(q-1)\zeta_i\equiv0\pmod{1}$ for $i=1,2$, and
\[
\begin{array}{l}
\left(\begin{array}{cc}
m_1   & m_1-m \\
n_1-n & n_1
\end{array}\right)
\left(\begin{array}{c}
\zeta_1\\
\zeta_2
\end{array}\right)
\equiv0\pmod{1}.
\end{array}
\]
and
\[
\zeta_1\not\equiv0\pmod{1},\zeta_2\not\equiv0\pmod{1},\zeta_1+\zeta_2\not\equiv0\pmod{1}.
\]

Let $p=char(\mathbb{F}_{q}), \zeta=\frac{a}{b}\in\mathbb{Q}$ with $(a,b)=1$, note that $(q^r-1)\zeta\equiv0\pmod{1}$ for some positive integer $r$ if and only if $p\nmid b$. In other words, let $v_p(\cdot)$ be the $p$-adic valuation of $\mathbb{Q}$, then $(q^r-1)\zeta\equiv0\pmod{1}$ for some $r\in \mathbb{N}$ if and only if $v_p(\zeta)\geq 0$. On the other hand, let $m=p^km'\in \mathbb{Z}, (p,m')=1,k\geq0$, then $(m,q^r-1)=(m'(p^k,q^r-1),q^r-1)=(m',q^r-1)$. Hence we will denote $m'$ by $m_p$ in the following calculation of zeta function.

Now let $\bar{C}$ be the projective curve given by the equation of the affine curve $C$. Then the number $\bar{N}$ of rational points over $\mathbb{F}_{q}$, on the curve $\bar{C}$, is related to the number $N$ by $\bar{N}=N+2$, since $\bar{C}$ has 2 points at infinity, namely $[1:0:0]$ and $[0:1:0]$. We are going to calculate the zeta function of the curve $\bar{C}$. We will use the following notation as in \cite{Weil}:

\begin{notation}
$\bar{N}_v$ is the number of points of the curve $\bar{C}$ over the extension of degree $v$ of the ground-field $\mathbb{F}_{q}$. Let $\alpha_1,\alpha_2,\cdots,\alpha_r$ be a set of rational numbers and $\alpha_i\not\equiv0\pmod{1}$ for $i=1,,\cdots,r$. We denote by $\mu=\mu(\alpha)=\mu(\alpha_1,\alpha_2,\cdots,\alpha_r)$ be the smallest positive integer such that $(q^\mu-1)\alpha_i\equiv0\pmod{1}$ for $1\leq i \leq r$. Suppose $\chi_{\alpha_1},\chi_{\alpha_2},\cdots,\chi_{\alpha_r}$ be multiplicative characters in the extension field $k$ of degree $\mu(\alpha)$ of the ground field $\mathbb{F}_{q}$, we denote by \linebreak $\chi_{\alpha_{i,\lambda}},g(\psi,\chi_{\alpha_{i,\lambda}}),j(\chi_{\alpha_{1,\lambda}},\cdots,\chi_{\alpha_{r,\lambda}}),
j_0(\chi_{\alpha_{1,\lambda}},\cdots,\chi_{\alpha_{r,\lambda}})$ the corresponding \linebreak characters and sums for the extension of degree $\lambda$ of the field $k$.
\end{notation}

\begin{remark}
Let $\alpha_1,\alpha_2,\cdots,\alpha_r$ be a set of rational numbers and $\alpha_i\not\equiv0\linebreak\pmod{1}$ for $i=1,\cdots,r$. Let $\omega$ be a generator of $\mathbb{F}_{q^{\mu(\alpha)}}^*$. Let $\chi$ be a multiplicative character on $\mathbb{F}_{q^{\mu(\alpha)}}^*$ defined by $\chi(\omega)=e^{2\pi i\frac{1}{q^{\mu(\alpha)}-1}}$. Since $(q^{\mu(\alpha)}-1)\alpha_i\equiv0\pmod{1}$, we can find $r$ characters $\chi_1,\cdots,\chi_r$ on $\mathbb{F}_{q^{\mu(\alpha)}}^*$ such that $\chi_1(\omega)=e^{2\pi i\alpha_1},\cdots,\chi_r(\omega)=e^{2\pi i\alpha_r}$. Using former notation we can denote $\chi_1,\cdots,\chi_r$ by $\chi_{\alpha_1,\omega},\cdots,\chi_{\alpha_r,\omega}$ or $\chi_{\alpha_1},\cdots,\chi_{\alpha_r}$ to omit $\omega$ if not confused. Hence whenever we talk about $\chi_{\alpha_1},\cdots,\chi_{\alpha_r}$ on a finite field $\mathbb{F}_{q}^*$, we always mean that we have chosen a generator $\omega$ of $\mathbb{F}_{q}^*$ such that $\chi_{\alpha_i}(\omega)=e^{2\pi i\alpha_i}$ for $i=1,\cdots,r$.
\end{remark}

Thus
\[
\begin{array}{l}
\sum\limits_{v=1}^{\infty}\bar{N}_vU^{v-1}\\
=\sum\limits_{v=1}^{\infty}(q^v+1)U^{v-1}-\sum\limits_{v=1}^{\infty}
\sum\limits_{(m_1-m,n_1,q^v-1)\xi_{1,v}\equiv0\pmod{1}\atop \xi_{1,v}\not\equiv0\pmod{1}}\chi_{\xi_{1,v}}(-k_1^{-1})U^{v-1}\\
\quad-\sum\limits_{v=1}^{\infty}\sum\limits_{(m_1,n_1-n,q^v-1)\xi_{2,v}\equiv0\pmod{1}\atop\xi_{2,v}\not\equiv0\pmod{1}}\chi_{\xi_{2,v}}(-k_2^{-1})U^{v-1}\\
\quad-\sum\limits_{v=1}^{\infty}\sum\limits_{(m,n,q^v-1)\xi_{3,v}\equiv0\pmod{1}\atop\xi_{3,v}\not\equiv0\pmod{1}}\chi_{\xi_{3,v}}(-k_1k_2^{-1})U^{v-1}\\
\quad+\sum\limits_{v=1}^{\infty}\frac{1}{q^v}\sum\limits_{\zeta_{1,v}\atop\zeta_{2,v}}\chi_{\zeta_{1,v}}(k_2^{-1})\chi_{\zeta_{2,v}}(k_1^{-1})
g(\psi,\chi_{\zeta_{1,v}})g(\psi,\chi_{\zeta_{2,v}})\\
\hspace{20em}\cdot g(\psi,\chi_{-\zeta_{1,v}-\zeta_{2,v}})U^{v-1}\\
\hspace{2em}\text{where $\zeta_{1,v},\zeta_{2,v}$ satisfy}
\begin{array}{l}
\left(\begin{array}{cc}
m_1   & m_1-m \\
n_1-n & n_1
\end{array}\right)
\left(\begin{array}{c}
\zeta_{1,v}\\
\zeta_{2,v}
\end{array}\right)
\equiv0\pmod{1}
\end{array}\\
%\end{array}
%\]
%\[
%\begin{array}{l}
\hspace{2em}\text{and}
\quad\zeta_{i,v}\not\equiv0\pmod{1},(q^v-1)\zeta_{i,v}\equiv0\pmod{1}, \text{for i=1,2},\\
\hspace{2em}\text{and}\quad \zeta_{1,v}+\zeta_{2,v}\not\equiv0\pmod{1}.\\
%\end{array}
%\]
%\[
%\begin{array}{l}
=\frac{q}{1-qU}+\frac{1}{1-U}-\sum\limits_{(m_1-m,n_1)_p\xi_1\equiv0\pmod{1}\atop \xi_1\not\equiv0\pmod{1}}\sum\limits_{\lambda=1}^{\infty}\chi_{\xi_{1,\lambda}}(-k_1^{-1})U^{\lambda\mu(\xi_1)-1}\\
\quad-\sum\limits_{(m_1,n_1-n)_p\xi_2\equiv0\pmod{1}\atop\xi_2\not\equiv0\pmod{1}}\sum\limits_{\lambda=1}^{\infty}\chi_{\xi_{2,\lambda}}(-k_2^{-1})U^{\lambda\mu(\xi_2)-1}\\
\quad-\sum\limits_{(m,n)_p\xi_3\equiv0\pmod{1}\atop\xi_3\not\equiv0\pmod{1}}\sum\limits_{\lambda=1}^{\infty}\chi_{\xi_{3,\lambda}}(-k_1k_2^{-1})U^{\lambda\mu(\xi_3)-1}\\
\quad+\sum\limits_{\zeta_1,\zeta_2}\sum\limits_{\lambda=1}^{\infty}\frac{1}{q^{\lambda\mu(\zeta_1,\zeta_2)}}\chi_{\zeta_{1,\lambda}}(k_2^{-1})\chi_{\zeta_{2,\lambda}}(k_1^{-1})\\
\hspace{6em}\cdot g(\psi,\chi_{\zeta_{1,\lambda}})g(\psi,\chi_{\zeta_{2,\lambda}})g(\psi,\chi_{-\zeta_{1,\lambda}-\zeta_{2,\lambda}})U^{\lambda\mu(\zeta_1,\zeta_2)-1}\\
\linebreak\\
\hspace{2em}\text{where $\zeta_1,\zeta_2$ satisfy}
\begin{array}{l}
\left(\begin{array}{cc}
m_1   & m_1-m \\
n_1-n & n_1
\end{array}\right)
\left(\begin{array}{c}
\zeta_1\\
\zeta_2
\end{array}\right)
\equiv0\pmod{1}
\end{array}
\\
\hspace{2em}\text{and}
\quad\zeta_i\not\equiv0\pmod{1},v_p(\zeta_i)\geq0,\text{ for i=1,2, } \zeta_1+\zeta_2\not\equiv0\pmod{1}\\
\end{array}
\]
\[
\begin{array}{l}
=\frac{q}{1-qU}+\frac{1}{1-U}-\sum\limits_{\xi_1}\sum\limits_{\lambda=1}^{\infty}\chi_{\xi_1}(-k_1^{-1})^\lambda U^{\lambda\mu(\xi_1)-1}\\
\quad-\sum\limits_{\xi_2}\sum\limits_{\lambda=1}^{\infty}\chi_{\xi_2}(-k_2^{-1})^\lambda U^{\lambda\mu(\xi_2)-1}
-\sum\limits_{\xi_3}\sum\limits_{\lambda=1}^{\infty}\chi_{\xi_3}(-k_1k_2^{-1})^\lambda U^{\lambda\mu(\xi_3)-1}\\
\quad+\sum\limits_{\zeta_1,\zeta_2}\sum\limits_{\lambda=1}^{\infty}\frac{1}{q^{\lambda\mu(\zeta_1,\zeta_2)}}\chi_{\zeta_1}(k_2^{-1})^\lambda\chi_{\zeta_2}(k_1^{-1})^\lambda\\
\hspace{11em}\cdot g(\psi,\chi_{\zeta_1})^\lambda g(\psi,\chi_{\zeta_2})^\lambda g(\psi,\chi_{-\zeta_1-\zeta_2})^\lambda U^{\lambda\mu(\zeta_1,\zeta_2)-1}\\
=\frac{q}{1-qU}+\frac{1}{1-U}-\sum\limits_{\xi_1}\frac{\chi_{\xi_1}(-k_1^{-1})U^{\mu(\xi_1)-1}}{1-\chi_{\xi_1}(-k_1^{-1})U^{\mu(\xi_1)}}
-\sum\limits_{\xi_2}\frac{\chi_{\xi_2}(-k_2^{-1})U^{\mu(\xi_2)-1}}{1-\chi_{\xi_2}(-k_2^{-1})U^{\mu(\xi_2)}}\\
\quad-\sum\limits_{\xi_3}\frac{\chi_{\xi_3}(-k_1k_2^{-1})U^{\mu(\xi_3)-1}}{1-\chi_{\xi_3}(-k_1k_2^{-1})U^{\mu(\xi_3)}}
-\sum\limits_{\zeta_1,\zeta_2}\frac{C(\zeta_1,\zeta_2)U^{\mu(\zeta_1,\zeta_2)-1}}{1-C(\zeta_1,\zeta_2)U^{\mu(\zeta_1,\zeta_2)}}\\
\hspace{0.5em}\text{where $C(\zeta_1,\zeta_2)=-\frac{1}{q^{\mu(\zeta_1,\zeta_2)}}\chi_{\zeta_1}(k_2^{-1})\chi_{\zeta_2}(k_1^{-1})
g(\psi,\chi_{\zeta_1})g(\psi,\chi_{\zeta_2})$}\\
\hspace{20em}\cdot g(\psi,\chi_{-\zeta_1-\zeta_2})\\
%\linebreak\\
\end{array}
\]
\[
\begin{array}{l}
=-\frac{d}{dU}\log(1-qU)-\frac{d}{dU}\log(1-U)+\sum\limits_{\xi_1}\frac{1}{\mu(\xi_1)}\frac{d}{dU}\log(1-\chi_{\xi_1}(-k_1^{-1})U^{\mu(\xi_1)})\\
\quad+\sum\limits_{\xi_2}\frac{1}{\mu(\xi_2)}\frac{d}{dU}\log(1-\chi_{\xi_2}(-k_2^{-1})U^{\mu(\xi_2)})\\
\quad+\sum\limits_{\xi_3}\frac{1}{\mu(\xi_3)}\frac{d}{dU}\log(1-\chi_{\xi_3}(-k_1k_2^{-1})U^{\mu(\xi_3)})\\
\quad+\sum\limits_{\zeta_1,\zeta_2}\frac{1}{\mu(\zeta_1,\zeta_2)}\frac{d}{dU}\log(1-C(\zeta_1,\zeta_2)U^{\mu(\zeta_1,\zeta_2)})\\
=\frac{d}{dU}\log(\frac{1}{(1-U)(1-qU)}\cdot\prod\limits_{(m_1-m,n_1)_p\xi_1\equiv0\pmod{1}\atop \xi_1\not\equiv0\pmod{1}}(1-\chi_{\xi_1}(-k_1^{-1})U^{\mu(\xi_1)})\\
\quad\cdot\prod\limits_{(m_1,n_1-n)_p\xi_2\equiv0\pmod{1}\atop\xi_2\not\equiv0\pmod{1}}(1-\chi_{\xi_2}(-k_2^{-1})U^{\mu(\xi_2)})\\
\quad\cdot\prod\limits_{(m,n)_p\xi_3\equiv0\pmod{1}\atop\xi_3\not\equiv0\pmod{1}}(1-\chi_{\xi_3}(-k_1k_2^{-1})U^{\mu(\xi_3)})\\
\quad\cdot\prod\limits_{\zeta_1,\zeta_2}(1-C(\zeta_1,\zeta_2)U^{\mu(\zeta_1,\zeta_2)})
)\\
\hspace{2em}\text{where $\zeta_1,\zeta_2$ satisfy}
\begin{array}{l}
\left(\begin{array}{cc}
m_1   & m_1-m \\
n_1-n & n_1
\end{array}\right)
\left(\begin{array}{c}
\zeta_1\\
\zeta_2
\end{array}\right)
\equiv0\pmod{1}
\end{array}
\\
\hspace{2em}\text{and
$\quad\zeta_i\not\equiv0\pmod{1},v_p(\zeta_i)\geq0,\text{ for i=1,2, } \zeta_1+\zeta_2\not\equiv0\pmod{1}$ but}\\
\hspace{2em}\text{taking only one representative for each set of pairs $(q^\rho\zeta_1,q^\rho\zeta_2)$ with}\\
\hspace{2em}\text{$0\leq \rho < \mu(\zeta_1,\zeta_2)$. Similarly for $\xi_i$}.\\

\end{array}
\]

Hence by \cite[Theorem 2.1]{Aubry}, the numerator of the zeta function of the nonsingular model $\tilde{C}$ of the affine curve $C$ is
\[
\begin{array}{l}
P_{\tilde{C}}(U)=\prod\limits_{\zeta_1,\zeta_2}(1-C(\zeta_1,\zeta_2)U^{\mu(\zeta_1,\zeta_2)})\\
=\prod\limits_{\zeta_1,\zeta_2}(1+\frac{1}{q^{\mu(\zeta)}}\chi_{\zeta_1}(k_2^{-1})\chi_{\zeta_2}(k_1^{-1})
g(\psi,\chi_{\zeta_1})g(\psi,\chi_{\zeta_2})g(\psi,\chi_{-\zeta_1-\zeta_2})U^{\mu(\zeta)})
\end{array}
\]
the product being taking over all $\zeta=(\zeta_1,\zeta_2)^t$ satisfying
\[
\begin{array}{l}
\left(\begin{array}{cc}
m_1   & m_1-m \\
n_1-n & n_1
\end{array}\right)
\left(\begin{array}{c}
\zeta_1\\
\zeta_2
\end{array}\right)
\equiv0\pmod{1}
\end{array}
\]
and $\zeta_i\not\equiv0\pmod{1},v_p(\zeta_i)\geq0,\text{ for i=1,2, } \zeta_1+\zeta_2\not\equiv0\pmod{1}$ but taking only one representative for each set of pairs $(q^\rho\zeta_1,q^\rho\zeta_2)$ with $0\leq \rho < \mu(\zeta_1,\zeta_2)$.

\begin{lemma}\label{lemma:mu}
Let $\alpha=(\alpha_1,\alpha_2,\cdots,\alpha_n)^t,\beta=(\beta_1,\beta_2,\cdots,\beta_n)^t$ with $\alpha_i,\beta_i \in\mathbb{Q}$ for $1\leq i \leq n$. Let $p$ be a prime number, denote $v_p(\alpha_i)\ge 0$ for all $1\le i\le n$ by $v_p(\alpha)\ge 0$. Let $m_\alpha$ be the smallest positive integer such that $m_\alpha\alpha_i\equiv0\pmod{1}$ for $i=1,2,\cdots,n$. Similarly for $m_\beta$. Suppose $\alpha=U\beta$ with $U=(u_{i,j})_{n\times n}\in GL(n,\mathbb{Z})$, then $\mu(\alpha)=\mu(\beta)$, $m_\alpha=m_\beta$, and $v_p(\alpha)\ge 0$ if and only if $v_p(\beta)\ge 0$.
\end{lemma}
\proof
Suppose $(q^{\mu(\beta)}-1)\beta_i\equiv0\pmod{1}$ for all $1\leq i \leq n$, then $(q^{\mu(\beta)}-1)\alpha_i\equiv(q^{\mu(\beta)}-1)(\sum\limits_{j=1}^nu_{ij}\beta_j)\equiv0\pmod{1}$ for all $1\leq i \leq n$. Hence $\mu(\alpha)\leq\mu(\beta)$. Since $\beta=U^{-1}\alpha$, we also have $\mu(\beta)\leq\mu(\alpha)$. Therefore $\mu(\alpha)=\mu(\beta)$. Similarly we have $m_\alpha=m_\beta$.

Let $\alpha_i=\frac{r_i}{s_i},(r_i,s_i)=1$, then $m_\alpha$ is the least common multiple of\linebreak $s_1,s_2,\cdots,s_n$. Since $v_p(\alpha)\ge 0\Leftrightarrow p\nmid m_\alpha$ and $m_\alpha=m_\beta$, we have $v_p(\alpha)\ge 0$ if and only if $v_p(\beta)\ge 0$.
\qed

\begin{prop}\label{prop:genus-m(C)}
Let $A=\left(\begin{array}{cc}
a_{11}  &a_{12}\\
a_{21}  &a_{22}
\end{array}\right)$
be a $2\times2$ matrix with $a_{ij}\in\mathbb{Z}$ and the determinant $|A|\ne0$. Let $d=(a_{11},a_{12},a_{21},a_{22})$. Let $C$ be a projective, non-singular, geometrically irreducible algebraic curve defined over a finite field $\mathbb{F}_q$, and let $p=char(\mathbb{F}_{q})$. Suppose the numerator of the zeta function of the curve $C$ is of the form
\[
\begin{array}{l}
P_{C}(U)
=\prod\limits_{\zeta_1,\zeta_2}(1+k(\zeta_1,\zeta_2)U^{\mu(\zeta_1,\zeta_2)})
\end{array}
\]
the product being taking over all $\zeta=(\zeta_1,\zeta_2)^t$ satisfying $A\zeta\equiv0\pmod{1}$
and $\zeta_i\not\equiv0\pmod{1},v_p(\zeta_i)\geq0,\text{ for i=1,2, } \zeta_1+\zeta_2\not\equiv0\pmod{1}$ but taking only one representative for each set of pairs $(q^\rho\zeta_1,q^\rho\zeta_2)$ with $0\leq \rho < \mu(\zeta_1,\zeta_2)$. And for any such pair $(\zeta_1,\zeta_2)$, $k(\zeta_1,\zeta_2)$ is a nonzero complex number. Let $m_C:=\max\{m_\zeta:\text{for all pairs $\zeta=(\zeta_1,\zeta_2)^t$ in the product of $P_{C}(U)$}\}$. Then \begin{enumerate}
\item[(1)] the genus $g(C)$ of $C$ over $\mathbb{F}_q$ is
\[
\frac{||A||_p-(a_{11},a_{21})_p-(a_{12},a_{22})_p-(a_{11}-a_{12},a_{21}-a_{22})_p}{2}+1,
\]
where $||A||$ is the absolute value of $|A|$.

\item[(2)] Suppose $p\nmid d$ and the genus $g(C)>0$, then $m_C=(\frac{||A||}{d})_p$.
\end{enumerate}

\end{prop}

\proof (1) We know that (\cite{Smith} or \cite[Chapter 14]{Hua2}) the Smith normal form of the matrix $A$ is
$
B=\left(\begin{array}{cc}
d  &0 \\
0  &\frac{|A|}{d}
\end{array}\right)
$, that's to say, we can find matrices $U,V\in GL(2,\mathbb{Z})$ such that
$UAV=B$. Let $V^{-1}\zeta=\eta=(\eta_1,\eta_2)^t$, then $d\eta_1\equiv0\pmod{1},\frac{|A|}{d}\eta_2\equiv0\pmod{1}$.
Since $\mu(\zeta)=\mu(V^{-1}\zeta)=\mu(\eta)$ by Lemma \ref{lemma:mu}, we have
\[
\sum_{A\zeta\equiv0\pmod{1}\atop v_p(\zeta_i)\geq0,\text{ for i=1,2, }}\mu(\zeta)=
\sum_{B\eta\equiv0\pmod{1}\atop v_p(\eta_i)\geq0,\text{ for i=1,2, }}\mu(\eta)
=d_p\cdot|\frac{|A|}{d}|_p=||A||_p,
\]
where $\zeta$'s in the sum take only one representative for each set of pairs $(q^\rho\zeta_1,q^\rho\zeta_2)$ with $0\leq \rho < \mu(\zeta_1,\zeta_2)$, similarly for $\eta$.
Next we count the sum of $\mu(\zeta)$ for $\zeta$ satisfying $A\zeta\equiv0\pmod{1}$
and $\zeta_1\not\equiv0\pmod{1}, \zeta_2\equiv0\pmod{1},v_p(\zeta_i)\geq0,\text{ for i=1,2, }$ but taking only one representative for each set of pairs $(q^\rho\zeta_1,q^\rho\zeta_2)$ with $0\leq \rho < \mu(\zeta_1,\zeta_2)$. Since $A\zeta\equiv0\pmod{1},\zeta_1\not\equiv0\pmod{1}, \zeta_2\equiv0\pmod{1},v_p(\zeta_i)\geq0,\text{ for i=1,2 }$ are equivalent to
$(a_{11},a_{21})_p\zeta_1\equiv0\pmod{1},\zeta_1\not\equiv0\pmod{1}, \zeta_2\equiv0\pmod{1}$, we see that the sum of $\mu(\zeta)$ for this case is\linebreak $(a_{11},a_{21})_p-1$. Similarly, the sum of $\mu(\zeta)$ for the case of $\zeta_1\equiv0\pmod{1}, \zeta_2\not\equiv0\pmod{1}$ is $(a_{12},a_{22})_p-1$, and for the case of $\zeta_1\not\equiv0\pmod{1}, \zeta_1+\zeta_2\equiv0\pmod{1}$ is $(a_{11}-a_{12},a_{21}-a_{22})_p-1$. Therefore,
$\deg P_C(U)=||A||_p-1-((a_{11},a_{21})_p-1)-((a_{12},a_{22})_p-1)-((a_{11}-a_{12},a_{21}-a_{22})_p-1)
$
and
\[
\begin{array}{ll}
g(C)&=\frac{\deg P_{C}(U)}{2}\\
&=\frac{||A||_p-(a_{11},a_{21})_p-(a_{12},a_{22})_p-(a_{11}-a_{12},a_{21}-a_{22})_p}{2}+1
\end{array}
\]
by \cite[Theorem 5.1.15]{Stichtenoth}.

(2) Let $N$ be the number of pairs $\eta=(\eta_1,\eta_2)^t$ satisfying $d\eta_1\equiv0\pmod{1},\linebreak \frac{|A|}{d}\eta_2\equiv0\pmod{1},v_p(\eta_i)\geq0,\text{ for i=1,2, }m_{\eta}=(\frac{||A||}{d})_p$, then $N\ge d\cdot\varphi((\frac{||A||}{d})_p)$ since $d|\frac{|A|}{d}$, where $\varphi$ is the Euler function. Thus in the case $d>3$ or in the case $\max\{(a_{11},a_{21})_p,(a_{12},a_{22})_p,(a_{11}-a_{12},a_{21}-a_{22})_p\}<(\frac{||A||}{d})_p$,  we can see that $m_C=(\frac{||A||}{d})_p$ by Lemma \ref{lemma:mu}. Now suppose $d\le 3$ and $\max\{(a_{11},a_{21})_p,(a_{12},\linebreak a_{22})_p, (a_{11}-a_{12},a_{21}-a_{22})_p\}=(\frac{||A||}{d})_p$. Assume that, for example, $d=2,(a_{11},a_{21})_p=(\frac{||A||}{d})_p,\max\{(a_{12},a_{22})_p, (a_{11}-a_{12},a_{21}-a_{22})_p\}<(\frac{||A||}{d})_p$. Note that the number of pairs $\zeta=(\zeta_1,\zeta_2)^t$ satisfying $(a_{11},a_{21})_p\zeta_1\equiv0\pmod{1},\zeta_1\not\equiv0\pmod{1}, \zeta_2\equiv0\pmod{1},v_p(\zeta_i)\geq0,\text{ for i=1,2, }m_{\zeta}=(\frac{||A||}{d})_p$ is $\varphi((\frac{||A||}{d})_p)$, we still see that $m_C=(\frac{||A||}{d})_p$ since $N\ge 2\cdot\varphi((\frac{||A||}{d})_p)$. So it is reduced to consider the case $d=3,(a_{11},a_{21})_p=(a_{12},a_{22})_p=(a_{11}-a_{12},a_{21}-a_{22})_p=(\frac{||A||}{d})_p$,  since for any other cases either $g(C)\le 0$ or we can easily see that $m_C=(\frac{||A||}{d})_p$.
For this case, we will show that $d=(\frac{||A||}{d})_p=3$. In fact, let $d_1=(a_{11},a_{21})=p^{j}m,\frac{||A||}{d}=\frac{|a_{11}a_{22}-a_{21}a_{12}|}{3}=p^im,p\nmid m$, then $j\le i$ since $d_1|\frac{||A||}{d}$, and we have $\frac{a_{11}}{d_1}\frac{a_{22}}{3}-\frac{a_{21}}{d_1}\frac{a_{12}}{3}=\pm p^{i-j}$. It follows that $(\frac{a_{12}}{3},\frac{a_{22}}{3})_p=1$, i.e. $(a_{12},a_{22})_p=3$. Thus the number of pairs of $\zeta=(\zeta_1,\zeta_2)^t$ satisfying $A\zeta\equiv0\pmod{1}$
and $\zeta_i\not\equiv0\pmod{1},v_p(\zeta_i)\geq0,\text{ for i=1,2, } \zeta_1+\zeta_2\not\equiv0\pmod{1}, m_{\zeta}=(\frac{||A||}{d})_p$ is at least $2$. Hence we also have $m_C=(\frac{||A||}{d})_p$ for this last case.
\qed

\begin{corollary}
Let $p=char(\mathbb{F}_{q})$, and let $C$ be a the nonsingular model over $\mathbb{F}_{q}$ of the geometrically irreducible curve with equation $x^{m_1}y^{n_1}+k_1x^m+k_2y^n=0$, where $m_1+n_1>m, m_1+n_1>n$, $n_1 \geq m_1$, and if $m_1=n_1$ then $n\geq m$. Let $g(C)$ be the genus of $C$. If $i(C)=0$, then $g(C)=0$. If $i(C)>0$, then the genus $g(C)$ of $C$ is
\[
\frac{(m_1n+mn_1-mn)_p-(m_1-m,n_1)_p-(m_1,n_1-n)_p-(m,n)_p}{2}+1.
\]
\end{corollary}
\proof
If $i(C)=0$, then $g(C)=0$ since $g(C)\le i(C)$ by Theorem \ref{thm:Baker}. If $i(C)>0$, from Proposition \ref{prop:i(C)} we know that $m_1n+mn_1-mn\ge0$.
\qed

Next we will show that the form of the numerator of zeta function before Lemma \ref{lemma:mu} has a similar form as in \cite{Tate}.

\begin{lemma}\cite{Tate}\label{lemma:Tate}
Let $k$ be a quadratic extension of the finite field $k_0$ of $q$ elements. Denote by $\theta$ a nontrivial multiplicative character of $k^*$ which is trivial on $k_0^*$ and by $\psi$ the standard additive character of $k$. Then
\[
g(\psi,\theta)=\theta(c)q
\]
for some $c\in k^*$ which satisfies $Tr_{k/k_0}(c)=0$.
\end{lemma}

\begin{prop}\label{prop:Tate_Shafa}
Let $A=\left(\begin{array}{cc}
a_{11}  &a_{12}\\
a_{21}  &a_{22}
\end{array}\right)$
be a $2\times2$ matrix with $a_{ij}\in\mathbb{Z}$ and the determinant $|A|\ne0$. Let $\xi=(\xi_1,\xi_2)$ be a pair of rational numbers such that
\[
\begin{array}{l}
\left\{\begin{array}{llll}
a_{11}\xi_1&+&a_{12}\xi_2&\equiv0\pmod{1}\\
a_{21}\xi_1&+&a_{22}\xi_2&\equiv0\pmod{1}
\end{array}\right.
\end{array}
\]
and $\xi_i\not\equiv0\pmod{1},v_p(\xi_i)\geq0$, for $i=1,2,$ $(\xi_1+\xi_2)\not\equiv0\pmod{1}$. \textup{\quad$(**)$}
Let $\omega$ be a generator of $\mathbb{F}_{q^{\mu(\xi)}}^*$ such that $\chi_{\xi_1},\chi_{\xi_2},\chi_{-\xi_1-\xi_2}$ are multiplicative characters on $\mathbb{F}_{q^{\mu(\xi)}}^*$ and $\chi_{\xi_i}(\omega)=e^{2\pi i\xi_i}$ for $i=1,2$. Let $k_1,k_2\in\mathbb{F}_{q}^*$. Let $d=(a_{11},a_{12},a_{21},a_{22})$. Let $p=char(\mathbb{F}_{q})$. Suppose $p\nmid d$ and $(\frac{|A|}{d})_p|(q^n+1)$ for some n. Then $\mu(\xi)$ is even and $\mu(\xi)=2(n,\mu(\xi))$. Let $\mu(\xi)=2\nu(\xi)$. Then
\[
\frac{1}{q^{\mu(\xi)}}\chi_{\xi_1}(k_1)\chi_{\xi_2}(k_2)
g(\psi,\chi_{\xi_1})g(\psi,\chi_{\xi_2})g(\psi,\chi_{-\xi_1-\xi_2})=q^{\nu(\xi)}.
\]
\end{prop}
\proof
The proof is similar to that of Theorem 1 in \cite{Tate}. We reproduce it here for the convenience of the reader.

We know that (\cite{Smith} or \cite[Chapter $14$]{Hua2}) the Smith normal form of the $2\times2$ matrix $A=\left(\begin{array}{cc}
a_{11}  &a_{12}\\
a_{21}  &a_{22}
\end{array}\right)$
is
$
\left(\begin{array}{cc}
d  &0 \\
0  &\frac{|A|}{d}
\end{array}\right)
$, that's to say, we can find matrices $U,V\in GL(2,\mathbb{Z})$ such that
$UAV=\left(\begin{array}{cc}
d  &0  \\
0  &\frac{|A|}{d}
\end{array}\right)$.
Let $V^{-1}\xi^t=\eta=(\eta_1,\eta_2)^t$, then $d\eta_1\equiv0\pmod{1},\frac{|A|}{d}\eta_2\equiv0\pmod{1}$. Since $d|\frac{|A|}{d}$, we have $\frac{|A|}{d}\xi_i\equiv0\pmod{1}$ for $i=1,2$. Since $v_p(\xi_i)\geq0$, then $(\frac{|A|}{d})_p\xi_i\equiv0\pmod{1}$ for $i=1,2$. Now given a pair of rational numbers $(\xi_1,\xi_2)$ satisfies the condition $(**)$, we first show that $\mu(\xi)$ is even.

Let $m_\xi$ be the smallest positive integer such that $m_\xi\xi_i\equiv0\pmod{1}$ for $i=1,2$. Then $m_\xi|(\frac{|A|}{d})_p$ and $\mu(\xi)$ is the smallest positive integer such that $m_\xi|(q^{\mu(\xi)}-1)$ by the definition of $\mu(\xi)$. Hence $q+m_\xi\mathbb{Z}$ has order $\mu(\xi)$ in $(\mathbb{Z}/m_\xi\mathbb{Z})^*$ and $q^n+m_\xi\mathbb{Z}$ has order $2$ since $m_\xi|(\frac{|A|}{d})_p$ and $(\frac{|A|}{d})_p|(q^n+1)$ by hypothesis, and $m_\xi>2$ by $(**)$. Let $\varphi$ be the Euler function, and let $(\mathbb{Z}/m_\xi\mathbb{Z})^*=<\tau>$. Since $(q,(\frac{|A|}{d})_p)=1$, it follows that $(q,m_\xi)=1$, so we can suppose that $q=\tau^a$ with $a\in\mathbb{Z}$. Thus $\mu(\xi)=\frac{\varphi(m_\xi)}{(a,\varphi(m_\xi))}$ and $2=\frac{\varphi(m_\xi)}{(na,\varphi(m_\xi))}$, which implies that $\mu(\xi)=2(n,\mu(\xi))$ and, in particular, that $\mu(\xi)$ is even.

Now let $G$ be the group of multiplicative characters on $\mathbb{F}_{q^{\mu(\xi)}}^*$. Let $H$ be the subgroup generated by $\chi_{\xi_1},\chi_{\xi_2}$ in $G$. Then the order $|H|$ of $H$ is the l.c.m. of the orders of $\chi_{\xi_1},\chi_{\xi_2}$ in $G$. Hence $|H|=m_\xi$. Note that $\nu(\xi)=\frac{1}{2}\mu(\xi)$. Next we will show that $\chi_\beta=1$ on $\mathbb{F}_{q^{\nu(\xi)}}^*$ for all $\chi_\beta\in H$. It suffices to show that $\chi_{\frac{1}{m_\xi}}=1$ on $\mathbb{F}_{q^{\nu(\xi)}}^*$. Note that
\[
\mathbb{F}_{q^{\nu(\xi)}}^*=<\omega^{\frac{q^{\mu(\xi)}-1}{q^{\nu(\xi)}-1}}>=<\omega^{q^{\nu(\xi)}+1}>
\quad\text{and}\quad
\chi_{\frac{1}{m_\xi}}(\omega^{q^{\nu(\xi)}+1})=e^{2\pi i\frac{q^{\nu(\xi)}+1}{m_\xi}}.
\]
Since $\nu(\xi)=(n,\mu(\xi))$, then $\frac{n}{\nu(\xi)}$ is odd and there exist integers $a,b$ with $a$ odd such that $\nu(\xi)=an+b\mu(\xi)$. Hence from $q^n\equiv -1\pmod{m_\xi}$ and $q^{\mu(\xi)}\equiv 1 \pmod{m_\xi}$, we have $q^{\nu(\xi)}\equiv -1 \pmod{m_\xi}$, i.e. $m_\xi|(q^{\nu(\xi)}+1)$. Thus $\chi_{\frac{1}{m_\xi}}(\omega^{q^{\nu(\xi)}+1})=1$.

Finally, for a character $\chi_\gamma$ on $\mathbb{F}_{q^{\mu(\xi)}}^*$, if $\chi_\gamma(\omega)=e^{2\pi i\gamma}$ and $\gamma \not\equiv0\pmod{1}$, then $\chi_\gamma$ is a nontrivial multiplicative character on $\mathbb{F}_{q^{\mu(\xi)}}^*$. It follows that $\chi_{\xi_1},\chi_{\xi_2}$ and $\chi_{-\xi_1-\xi_2}$ are all nontrivial on $\mathbb{F}_{q^{\mu(\xi)}}^*$. Therefore, from the lemma \ref{lemma:Tate} we see that
\[
\frac{1}{q^{\mu(\xi)}}\chi_{\xi_1}(k_1)\chi_{\xi_2}(k_2)
g(\psi,\chi_{\xi_1})g(\psi,\chi_{\xi_2})g(\psi,\chi_{-\xi_1-\xi_2})=q^{\nu(\xi)}.
\]
\qed

Summarizing the above results we have the following theorem.

\begin{theorem}\label{thm:zeta_function-case5}
Let $C$ be the nonsingular model over $\mathbb{F}_{q}$ of the geometrically irreducible curve given by
\[
x^{m_1}y^{n_1}+k_1x^m+k_2y^n=0,
\]
where $k_1,k_2\in\mathbb{F}_{q}^*; m_1+n_1>m, m_1+n_1>n$, and $n_1 \geq m_1$, if $m_1=n_1$ then $n\geq m$. Let $d=(m,n,m_1,n_1)$. Let $p=char(\mathbb{F}_{q})$. Suppose $p\nmid d$. Let $\xi=(\xi_1,\xi_2)$ be a pair of rational numbers such that
\[
\left\{\begin{array}{rcrl}
m_1\xi_1&+&(m_1-m)\xi_2&\equiv0\pmod{1}\\
(n_1-n)\xi_1&+&n_1\xi_2&\equiv0\pmod{1}
\end{array}\right.
\]
and
\[
\xi_i\not\equiv0\pmod{1},v_p(\xi_i)\geq0,\text{ for i=1,2, },
(\xi_1+\xi_2)\not\equiv0\pmod{1}.\textup{\quad$(***)$}
\]
Then
\begin{enumerate}
\item[(1)] The numerator of the zeta function of the curve $C$ is
\[
\begin{array}{l}
P_{\tilde{C}}(U)\\
=\prod\limits_{\xi}(1+\frac{1}{q^{\mu(\xi)}}\chi_{\xi_1}(k_2^{-1})\chi_{\xi_2}(k_1^{-1})
g(\psi,\chi_{\xi_1})g(\psi,\chi_{\xi_2})g(\psi,\chi_{-\xi_1-\xi_2})U^{\mu(\xi)}).
\end{array}
\]

\item[(2)] Suppose further that $((m_1n+m(n_1-n))/d)_p|(q^l+1)$ for some $l$. Then $\mu(\xi)$ is even and $\mu(\xi)=2(l,\mu(\xi))$. Let $\mu(\xi)=2\nu(\xi)$, then the numerator of the zeta function of the curve $C$ is
\[
P_{\tilde{C}}(U)=\prod_\xi(1+q^{\nu(\xi)}U^{\mu(\xi)}),
\]
\end{enumerate}
the product in (1) and (2) both being taking over all pairs $\xi=(\xi_1,\xi_2)$ satisfying $(***)$ but taking only one representative for each set of pairs $(q^\rho\xi_1,q^\rho\xi_2)$ with $0\leq \rho < \mu(\xi)$.

\end{theorem}

\qed

\begin{lemma}
Let $m$ be a positive integer and $q$ be a prime power. Suppose $(q,m)=1$ and the order of $q$ in $(\mathbb{Z}/m\mathbb{Z})^*$ is $k$, denoted by $|q|=k$, i.\ e.\ $k$ is the least positive integer such that $q^k\equiv 1\pmod{m}$. Then there exist $n$ such that $q^n \equiv -1\pmod{m}$ if and only if $k$ is even and $q^{\frac{k}{2}} \equiv -1\pmod{m}$.
\end{lemma}
\proof
It suffices to show that if there exist $n$ such that $q^n \equiv -1\pmod{m}$, then $k$ is even and $q^{\frac{k}{2}} \equiv -1\pmod{m}$.

Since $|q|=k$ and $q^{2n}\equiv 1\pmod{m}$, we have $k|2n$. If k is odd, then $(2,k)=1$, it follows that $k|n$. Hence $q^n\equiv 1\pmod{m}$, contrary to the hypothesis. Thus $k$ is even, and from $k|2n$ we have $\frac{k}{2}|n$. Now suppose $n=\frac{k}{2}(2i)=ki$, then $q^n=(q^k)^i\equiv 1\pmod{m}$, contrary to the hypothesis. Thus $n=\frac{k}{2}(2i+1)=ki+\frac{k}{2}$. Hence $ q^{\frac{k}{2}}\equiv (q^k)^iq^{\frac{k}{2}}\equiv q^n\equiv  -1\pmod{m}$.

\qed

\begin{definition}
Let $C$ be a projective, non-singular, geometrically irreducible algebraic curve defined over $\mathbb{F}_{q}$. We will call the numerator of the zeta function of $C$ the $L-$ polynomial of $C/\mathbb{F}_{q}$ as in \cite{Stichtenoth}, i.e. $L(t):=(1-t)(1-qt)Z(t)$, where $Z(t)$ is the zeta function of $C$.

\end{definition}

\begin{prop}\label{prop:maximal-L-Poly}
Let $C$ be a projective, non-singular, geometrically irreducible algebraic curve defined over $\mathbb{F}_{q}$. Suppose the $L-$ polynomial of $C/\mathbb{F}_{q}$ has the form:
\[
L_{C}(t)=\prod_{i=1}^r(1+q^{\nu_i}t^{\mu_i}),
\]
where $\nu_i$ are positive integer and $\mu_i=2\nu_i$ for all $1\leq i \leq r$. Let $n$ be a positive integer such that $\mu_i=2(n,\mu_i)$ for all $1\leq i \leq r$. Then the curve $C$ is maximal over the finite field $\mathbb{F}_{q^{2n}}$.
\end{prop}

\proof
Since the $L$-polynomial of $C$ over the finite field $\mathbb{F}_{q}$ is $L_C(t)=\prod_{i=1}^r(1+q^{\nu_i}t^{\mu_i})=\prod_{j=1}^{2g}(1-\alpha_jt)$, where g is the genus of the curve $C$, then the $L$-polynomial of $C$ over $\mathbb{F}_{q^{2n}}$ is $L_{2n}(t)=\prod_{j=1}^{2g}(1-\alpha_j^{2n}t)$. As we know, $C$ is maximal over $\mathbb{F}_{q^{2n}} \Leftrightarrow L_{2n}(t)=\prod_{j=1}^{2g}(1+q^nt) \Leftrightarrow -\alpha_j^{2n}=q^n$, i.\ e.\ $\alpha_j^{2n}=-q^n$ for all $j$, but $t=\frac{1}{\alpha_j}$ is a zero of $1+q^{\nu_i} t^{\mu_i}=0$ for some $\nu_i$. It follows that for each $j$, $\alpha_j^{2\nu_i}=-q^{\nu_i}$ for some $\nu_i$. Hence now it suffices to show that $\nu_i|n$ and $\frac{n}{\nu_i}$ are odd for all $\nu_i$. Note that $\mu_i=2(n,\mu_i)$ and $\mu_i=2\nu_i$, we have $\nu_i=(n,\mu_i)=\nu_i(\frac{n}{\nu_i},2)$. It follows that $\nu_i|n$ and $\frac{n}{\nu_i}$ are odd for all $\nu_i$ indeed.

\qed

\begin{theorem} \label{thm:maximal_case5}
Let $C$ be the nonsingular model of the geometrically irreducible curve defined over $\mathbb{F}_q$ by
\[
x^{m_1}y^{n_1}+k_1x^m+k_2y^n=0,
\]
where $k_1,k_2\in\mathbb{F}_q,k_1k_2\neq0; m_1+n_1>m, m_1+n_1>n$, and $n_1 \geq m_1$, if $m_1=n_1$ then $n\geq m$. Let $d=(m,n,m_1,n_1)$. Let $p=char(\mathbb{F}_{q})$. Suppose $p\nmid d$ and $((m_1n+m(n_1-n))/d)_p|(q+1)$. Then $C$ is maximal over $\mathbb{F}_{q^2}$. Conversely, if $C$ is maximal over $\mathbb{F}_{q^2}$ and $g(C)>0$, then $((m_1n+m(n_1-n))/d)_p|(q^2-1)$.
\end{theorem}

\proof Taking $n=1$, we can show that if $p\nmid d$ and $((m_1n+m(n_1-n))/d)_p|(q+1)$, then $C$ is maximal over $\mathbb{F}_{q^2}$ by Proposition \ref{prop:Tate_Shafa} and Proposition \ref{prop:maximal-L-Poly}. For the converse, if $C$ is maximal over $\mathbb{F}_{q^2}$ and $g(C)>0$, then $\mu(\xi)=1$ (cf. Theorem \ref{thm:zeta_function-case5} (1)) for all $\xi$ in the product of the numerator of the zeta function of the curve $C$ over $\mathbb{F}_{q^2}$. Since $m(C)=((m_1n+m(n_1-n))/d)_p$ by Proposition \ref{prop:genus-m(C)}, we have $((m_1n+m(n_1-n))/d)_p|(q^2-1)$.
\qed

\subsection{Other cases}
\subsubsection{Case 1}
For the Fermat curve with equation $k_1x^n+k_2y^n+1=0$ defined over $\mathbb{F}_{q}$, where $k_1,k_2\in\mathbb{F}_{q},k_1k_2\neq0,n\ge 2$, it is well-known \cite{Lachaud1} that if $n|(q+1)$ then $C$ is maximal over $\mathbb{F}_{q^2}$.

Let $C$ be the curve given by $k_1x^m+k_2y^n+1=0$ over $\mathbb{F}_{q}$, where $k_1,k_2\in\mathbb{Z},k_1k_2\neq0, n>m$. Let $d=(m,n),p=char(\mathbb{F}_{q})$. Suppose $p\nmid d$ and $(mn/d)|(q+1)$, we can show that $C$ is maximal over $\mathbb{F}_{q^2}$ by \cite[ Theorem 1]{Tate} and Proposition \ref{prop:maximal-L-Poly}. See also \cite[Theorem 5]{Tafazolian2}.

\subsubsection{Case 2}
Let $C$ be an geometrically irreducible affine curve over a finite field $\mathbb{F}_{q}$ with equation $k_1x^m+k_2y^{n_1}+y^{n}=0$, where $k_1k_2\in \mathbb{F}_{q}^*, n>m, n>n_1 $. Let $(n_1,n)=d, n_1s+nt=d,n_1=n_1'd,n=n'd,(m,q-1)=d_1,(d,q-1)=d_2,(n',q-1)=l_1,(n_1',q-1)=l_2$. Let $N$ be the number of points in $\mathbb{F}_{q}$ of the curve $C$. Put $L(u)=k_1u_0+k_2u_1+u_2$. Then
\[
\begin{array}{ll}
N&=\sum\limits_{L(u)=0}N_m(u_0)N_{n_1,n}(u_1,u_2)\\
&=\sum\limits_{L(u)=0}N_m(u_0)N_{d}(u_1^su_2^t)\delta(u_1^{n'},u_2^{n_1'})\\
&=1+\sum\limits_{L(u)=0\atop u_i\in\mathbb{F}_{q}^*}N_m(u_0)N_{d}(u_1^su_2^t)\delta(u_1^{n'},u_2^{n_1'})\\
&\quad+\sum\limits_{k_2u_1+u_2=0\atop u_1,u_2\in\mathbb{F}_{q}^*}N_{d}(u_1^su_2^t)\delta(u_1^{n'},u_2^{n_1'})\\
&=\cdots\\
&=1+\frac{1}{q-1}\sum\limits_{\theta_1,\theta_2}\chi_{\theta_1}(k_1^{-1})\chi_{\theta_2}(k_2^{-1})
j_0(\chi_{\theta_1},\chi_{\theta_2},\chi_{-\theta_1-\theta_2})\\
&\hspace{2em}+\frac{1}{q-1}\sum\limits_{(q-1,n-n_1)\theta_3\equiv0\pmod{1}}\chi_{\theta_3}(k_2^{-1})j_0(\chi_{\theta_3},\chi_{-\theta_3}),
\end{array}
\]
\hspace{3em} where $\theta_1,\theta_2$ in the first sum satisfy $(q-1)\theta_i\equiv0\pmod{1}$ for $i=1,2$ and
\[
\begin{array}{l}
\left(\begin{array}{cc}
m   &0 \\
n  &n-n_1
\end{array}\right)
\left(\begin{array}{c}
\theta_1\\
\theta_2
\end{array}\right)
\equiv0\pmod{1}.
\end{array}
\]

Thus similar to the case 5, we have

\begin{theorem}
Let $C$ be the nonsingular model over $\mathbb{F}_{q}$ of the geometrically irreducible curve given by
\[
k_1x^m+k_2y^{n_1}+y^{n}=0,
\]
where $k_1,k_2\in\mathbb{F}_{q}^*; n>m, n>n_1$. Let $d=(m,n,n_1)$. Let $p=char(\mathbb{F}_{q})$. Suppose $p\nmid d$. Let $i(C)=\frac{(m-1)(n-n_1)-(d_1+d_2)}{2}+1$, where $d_1=(m,n_1),d_2=(m,n)$. Let $g(C)$ be the genus of $C$ over $\mathbb{F}_{q}$. Let $\xi=(\xi_1,\xi_2)$ be a pair of rational numbers such that
\[
\left\{\begin{array}{rcrl}
m\xi_1&&&\equiv0\pmod{1}\\
n\xi_1&+&(n-n_1)\xi_2&\equiv0\pmod{1}
\end{array}\right.
\]
and
\[
\xi_i\not\equiv0\pmod{1},v_p(\xi_i)\geq0,\text{ for i=1,2, },
(\xi_1+\xi_2)\not\equiv0\pmod{1}.\textup{\quad$(4.1)$}
\]
Then
\begin{enumerate}
\item[(1)] If $i(C)=0$, then $g(C)=0$. If $i(C)>0$, then
\[
g(C)=\frac{(m(n-n_1))_p-(n-n_1)_p-(d_1)_p-(d_2)_p}{2}+1.
\]

\item[(2)] The numerator of the zeta function of the curve $C$ over $\mathbb{F}_{q}$ is
\[
\begin{array}{l}
P_{C}(U)\\
=\prod\limits_{\xi}(1+\frac{1}{q^{\mu(\xi)}}\chi_{\xi_1}(k_1^{-1})\chi_{\xi_2}(k_2^{-1})
g(\psi,\chi_{\xi_1})g(\psi,\chi_{\xi_2})g(\psi,\chi_{-\xi_1-\xi_2})U^{\mu(\xi)}).
\end{array}
\]

\item[(3)] Suppose further that $(m(n-n_1)/d)_p|(q^l+1)$ for some $l$. Then $\mu(\xi)$ is even and $\mu(\xi)=2(l,\mu(\xi))$. Let $\mu(\xi)=2\nu(\xi)$, then the numerator of the zeta function of the curve $C$ over $\mathbb{F}_{q}$ is
\[
P_{C}(U)=\prod_\xi(1+q^{\nu(\xi)}U^{\mu(\xi)}),
\]

the product in (2) and (3) both being taking over all pairs $\xi=(\xi_1,\xi_2)$ satisfying $(4.1)$ but taking only one representative for each set of pairs $(q^\rho\xi_1,q^\rho\xi_2)$ with $0\leq \rho < \mu(\xi)$.

\item[(4)] Suppose $((m(n-n_1)/d)_p|(q+1)$, i.e. $l=1$ in (3). Then $C$ is maximal over $\mathbb{F}_{q^2}$. Conversely, if $C$ is maximal over $\mathbb{F}_{q^2}$ and $g(C)>0$, then $((m(n-n_1)/d)_p|(q^2-1)$.

\end{enumerate}

\end{theorem}
\qed

\subsubsection{Case 3}
Let $C$ be an geometrically irreducible affine curve over a finite field $\mathbb{F}_{q}$ with equation $k_1x^{m_1}y^{n_1}+k_2y^{n}+1=0$, where $k_1k_2\in \mathbb{F}_{q}^*,n>m_1+n_1$. Let $(n_1,n)=d,n_1s+nt=d,n_1=n_1'd,n=n'd,(m_1,q-1)=d_1,(d,q-1)=d_2,(n_1',q-1)=l_1,(n',q-1)=l_2$. Let $N$ be the number of points in $\mathbb{F}_{q}$ of the curve $C$. Put $L(u,v)=k_1uv_0+k_2v_1+1$. Then
\[
\begin{array}{ll}
N&=\sum\limits_{L(u,v)=0}N_{m_1}(u)N_{n_1,n}(v_0,v_1)\\
&=\sum\limits_{L(u,v)=0}N_{m_1}(u)N_d(v_0^sv_1^t)\delta(v_0^{n'},v_1^{n_1'})\\
&=\sum\limits_{L(u,v)=0\atop u,v_i\in\mathbb{F}_{q}^*}N_{m_1}(u)N_{d}(v_0^sv_1^t)\delta(v_0^{n'},v_1^{n_1'})
+\sum\limits_{k_2v_1+1=0\atop v_0,v_1\in\mathbb{F}_{q}^*}N_d(v_0^sv_1^t)\delta(v_0^{n'},v_1^{n_1'})\\
&=\cdots\\
&=\sum\limits_{\theta_1,\theta_2}\chi_{\theta_1}(-k_1^{-1})\chi_{\theta_2}(-k_2^{-1})
j(\chi_{\theta_1},\chi_{\theta_2})+\sum\limits_{(q-1,n)\theta_3\equiv0\pmod{1}}\chi_{\theta_3}(-k_2^{-1})\\
&=\frac{1}{q-1}\sum\limits_{\theta_1,\theta_2}\chi_{\theta_1}(k_1^{-1})\chi_{\theta_2}(k_2^{-1})
j_0(\chi_{\theta_1},\chi_{\theta_2},\chi_{-\theta_1-\theta_2})\\
&\quad+\sum\limits_{(q-1,n)\theta_3\equiv0\pmod{1}}\chi_{\theta_3}(-k_2^{-1})\\
\end{array}
\]
\hspace{3em} where $\theta_1,\theta_2$ in the first sum satisfy $(q-1)\theta_i\equiv0\pmod{1}$ for $i=1,2$ and
\[
\begin{array}{l}
\left(\begin{array}{cc}
m_1   &0 \\
n_1  &n
\end{array}\right)
\left(\begin{array}{c}
\theta_1\\
\theta_2
\end{array}\right)
\equiv0\pmod{1}.
\end{array}
\]

Thus similarly we have

\begin{theorem}
Let $C$ be the nonsingular model over $\mathbb{F}_{q}$ of the geometrically irreducible curve given by
\[
k_1x^{m_1}y^{n_1}+k_2y^{n}+1=0,
\]
where $k_1,k_2\in\mathbb{F}_{q}^*; n>m_1+n_1$. Let $d=(m_1,n_1,n)$. Let $p=char(\mathbb{F}_{q})$. Suppose $p\nmid d$. Let $i(C)=\frac{(m_1-1)n-d_1-d_2}{2}+1$, where $d_1=gcd(m_1,n_1), d_2=gcd(m_1,n-n_1)$. Let $g(C)$ be the genus of $C$ over $\mathbb{F}_{q}$. Let $\xi=(\xi_1,\xi_2)$ be a pair of rational numbers such that
\[
\left\{\begin{array}{rcrl}
m_1\xi_1&&&\equiv0\pmod{1}\\
n_1\xi_1&+&n\xi_2&\equiv0\pmod{1}
\end{array}\right.
\]
and
\[
\xi_i\not\equiv0\pmod{1},v_p(\xi_i)\geq0,\text{ for i=1,2, },
(\xi_1+\xi_2)\not\equiv0\pmod{1}.\textup{\quad$(4.3)$}
\]
Then
\begin{enumerate}
\item[(1)] If $i(C)=0$, then $g(C)=0$. If $i(C)>0$, then
\[
g(C)=\frac{(m_1n)_p-n_p-(d_1)_p-(d_2)_p}{2}+1.
\]

\item[(2)] The numerator of the zeta function of the curve $C$ over $\mathbb{F}_{q}$ is
\[
\begin{array}{l}
P_{C}(U)\\
=\prod\limits_{\xi}(1+\frac{1}{q^{\mu(\xi)}}\chi_{\xi_1}(k_1^{-1})\chi_{\xi_2}(k_2^{-1})
g(\psi,\chi_{\xi_1})g(\psi,\chi_{\xi_2})g(\psi,\chi_{-\xi_1-\xi_2})U^{\mu(\xi)}).
\end{array}
\]

\item[(3)] Suppose further that $(m_1n/d)_p|(q^l+1)$ for some $l$. Then $\mu(\xi)$ is even and $\mu(\xi)=2(l,\mu(\xi))$. Let $\mu(\xi)=2\nu(\xi)$, then the numerator of the zeta function of the curve $C$ over $\mathbb{F}_{q}$ is
\[
P_{C}(U)=\prod_\xi(1+q^{\nu(\xi)}U^{\mu(\xi)}),
\]

the product in (2) and (3) both being taking over all pairs $\xi=(\xi_1,\xi_2)$ satisfying $(4.3)$ but taking only one representative for each set of pairs $(q^\rho\xi_1,q^\rho\xi_2)$ with $0\leq \rho < \mu(\xi)$.

\item[(4)] Suppose $(m_1n/d)_p|(q+1)$, i.e. $l=1$ in (3). Then $C$ is maximal over $\mathbb{F}_{q^2}$. Conversely, if $C$ is maximal over $\mathbb{F}_{q^2}$ and $g(C)>0$, then $(m_1n/d)_p|(q^2-1)$.

\end{enumerate}

\end{theorem}
\qed

\subsubsection{Case 4}
Let $C$ be an geometrically irreducible affine curve over a finite field $\mathbb{F}_{q}$ with equation $k_1x^{m_1}y^{n_1}+k_2x^{m}y^{n}+1=0$, where $k_1k_2\in \mathbb{F}_{q}^*,m_1+n_1\geq m+n, \frac{n_1}{m_1} \geq \frac{n}{m}$. Let $(m_1,m)=d_1,m_1=m_1'd_1,m=m'd_1,m_1s_1+mt_1=d_1,(n_1,n)=d_2,n_1=n_1'd_2,n=n'd_2,n_1s_2+nt_2=d_2,
(d_1,q-1)=d_1',(d_2,q-1)=d_2',(m',q-1)=l_1,(n',q-1)=l_2,(m_1',q-1)=l_3$ and $(n_1',q-1)=l_4$. Let $N$ be the number of points in $\mathbb{F}_{q}$ of the curve $C$. Put $L(u)=k_1u_0v_0+k_2u_1v_1+1$. Then
\[
\begin{array}{ll}
N&=\sum\limits_{L(u,v)=0}N_{m_1,m}(u_0,u_1)N_{n_1,n}(v_0,v_1)\\
&=\sum\limits_{L(u,v)=0}N_{d_1}(u_0^{s_1}u_1^{t_1})\delta(u_0^{m'},u_1^{m_1'})N_{d_2}(v_0^{s_2}v_1^{t_2})\delta(v_0^{n'},v_1^{n_1'})\\
&=\sum\limits_{L(u,v)=0\atop u_i,v_i\in\mathbb{F}_{q}^*}N_{d_1}(u_0^{s_1}u_1^{t_1})\delta(u_0^{m'},u_1^{m_1'})N_{d_2}(v_0^{s_2}v_1^{t_2})\delta(v_0^{n'},v_1^{n_1'})\\
&=\cdots\\
&=\sum\limits_{\theta_1,\theta_2}\chi_{\theta_1}(-k_2^{-1})\chi_{\theta_2}(-k_1^{-1})
j(\chi_{\theta_1},\chi_{\theta_2})\\
&=\frac{1}{q-1}\sum\limits_{\theta_1,\theta_2}\chi_{\theta_1}(k_2^{-1})\chi_{\theta_2}(k_1^{-1})
j_0(\chi_{\theta_1},\chi_{\theta_2},\chi_{-\theta_1-\theta_2})
\end{array}
\]
\hspace{3em} where $\theta_1,\theta_2$ satisfy $(q-1)\theta_i\equiv0\pmod{1}$ for $i=1,2$ and
\[
\begin{array}{l}
\left(\begin{array}{cc}
m   &m_1 \\
n  &n_1
\end{array}\right)
\left(\begin{array}{c}
\theta_1\\
\theta_2
\end{array}\right)
\equiv0\pmod{1}.
\end{array}
\]

Hence we have

\begin{theorem}
Let $C$ be the nonsingular model over $\mathbb{F}_{q}$ of the geometrically irreducible curve given by
\[
k_1x^{m_1}y^{n_1}+k_2x^{m}y^{n}+1=0,
\]
where $k_1,k_2\in\mathbb{F}_{q}^*;  m_1+n_1\geq m+n, \frac{n_1}{m_1} \geq \frac{n}{m}$. Let $d=(m,n,m_1,n_1)$. Let $p=char(\mathbb{F}_{q})$. Suppose $p\nmid d$. Let $i(C)=\frac{mn_1-m_1n-d_1-d_2-d_3}{2}+1$, where $d_1=gcd(m_1,n_1),d_2=gcd(m,n)$ and $d_3=gcd(n_1-n,m_1-m)$. Let $g(C)$ be the genus of $C$ over $\mathbb{F}_{q}$. Let $\xi=(\xi_1,\xi_2)$ be a pair of rational numbers such that
\[
\left\{\begin{array}{rcrl}
m\xi_1&+&m_1\xi_2&\equiv0\pmod{1}\\
n\xi_1&+&n_1\xi_2&\equiv0\pmod{1}
\end{array}\right.
\]
and
\[
\xi_i\not\equiv0\pmod{1},v_p(\xi_i)\geq0,\text{ for i=1,2, },
(\xi_1+\xi_2)\not\equiv0\pmod{1}.\textup{\quad$(4.3)$}
\]
Then
\begin{enumerate}
\item[(1)] If $i(C)=0$, then $g(C)=0$. If $i(C)>0$, then
\[
g(C)=\frac{(mn_1-m_1n)_p-(d_1)_p-(d_2)_p-(d_3)_p}{2}+1.
\]

\item[(2)] The numerator of the zeta function of the curve $C$ over $\mathbb{F}_{q}$ is
\[
\begin{array}{l}
P_{C}(U)\\
=\prod\limits_{\xi}(1+\frac{1}{q^{\mu(\xi)}}\chi_{\xi_1}(k_2^{-1})\chi_{\xi_2}(k_1^{-1})
g(\psi,\chi_{\xi_1})g(\psi,\chi_{\xi_2})g(\psi,\chi_{-\xi_1-\xi_2})U^{\mu(\xi)}).
\end{array}
\]

\item[(3)] Suppose further that $((mn_1-m_1n)/d)_p|(q^l+1)$ for some $l$. Then $\mu(\xi)$ is even and $\mu(\xi)=2(l,\mu(\xi))$. Let $\mu(\xi)=2\nu(\xi)$, then the numerator of the zeta function of the curve $C$ over $\mathbb{F}_{q}$ is
\[
P_{C}(U)=\prod_\xi(1+q^{\nu(\xi)}U^{\mu(\xi)}),
\]

the product in (2) and (3) both being taking over all pairs $\xi=(\xi_1,\xi_2)$ satisfying $(4.3)$ but taking only one representative for each set of pairs $(q^\rho\xi_1,q^\rho\xi_2)$ with $0\leq \rho < \mu(\xi)$.

\item[(4)] Suppose $(((mn_1-m_1n)/d)_p|(q+1)$, i.e. $l=1$ in (3). Then $C$ is maximal over $\mathbb{F}_{q^2}$. Conversely, if $C$ is maximal over $\mathbb{F}_{q^2}$ and $g(C)>0$, then $(((mn_1-m_1n)/d)_p|(q^2-1)$.

\end{enumerate}

\end{theorem}
\qed

\section*{Acknowledgements}

I started work on this paper when I was an exchange PhD student at the University of Groningen supported by the LiSUM Project in the framework of EMECW for 10 months in 2010. I wish to thank the project and the university's support and hospitality. I am particularly obliged to Prof. Jaap Top, my academic supervisor in Groningen. Particularly, this work was inspired by Prof. Jaap Top's example, $x^3+y^6+1=0$, which is of genus 4 and maximal over $\mathbb{F}_{5^2}$. In order to prove Theorem \ref{thm:maximal_case5}, Prof. Jaap Top suggested to study the zeta function and alerted me to \cite{Tate}.

\end{document}